\documentclass{amsart}

\usepackage{amsmath, amsthm, amssymb, amsfonts, amscd, verbatim, graphics}

\def\h{\hat{h}}

\def\0{{\bf 0}}
\def\B{{\bf B}}

\def\Sm{{\mathcal S}}

\def\f{{\mathfrak f}}

\def\F{{\mathbb F}}
\def\D{{\mathcal D}}
\def\Q{{\mathbb Q}}
\def\R{{\mathbb R}}
\def\Z{{\mathbb Z}}
\def\H{{\mathbb H}}
\def\C{{\mathbb C}}
\def\N{{\mathbb N}}
\def\H{{\mathbb H}}
\def\T{{\mathbb T}}
\def\O{{\mathcal O}}
\def\M{{\mathcal M}}

\def\Gal{\text{Gal}}

\def\tor{\text{tor}}
\def\ord{\text{ord}}
\def\tr{\text{tr}}

\theoremstyle{plain}

\newtheorem{thm}{Theorem}
\newtheorem{cor}[thm]{Corollary}
\newtheorem{prop}[thm]{Proposition}
\newtheorem{lem}[thm]{Lemma}

\def\hhat{\hat{h}}
\def\BB{{\bf B}}
\def\QQ{{\mathbb Q}}
\def\RR{{\mathbb R}}
\def\ZZ{{\mathbb Z}}

\def\CC{{\mathbb C}}

\def\G{{\mathbb G}}

\def\Qbar{\overline{\mathbb{Q}}}

\def\kbar{\overline{k}}
\def\Qpbar{\overline{\mathbb{Q}}_p}
\def\Fvbar{\overline{\mathbb{F}}_v}

\def\ee{{e}}
\def\ff{{f}}

\def\FF{{\mathbb F}}
\def\Ocal{{\mathcal O}}
\def\Ovhat{{\hat{\Ocal}_v}}

\newcommand{\Ecal}{{\mathcal E}}
\newcommand{\E}{{\mathcal E}}
\newcommand{\Ebar}{\bar{E}}
\newcommand{\Fpbar}{{\overline{\FF}_p}}
\newcommand{\Berk}{{\rm Berk}}
\newcommand{\PP}{\mathbb{P}}
\newcommand{\red}{{\pi}}

\newtheorem{theorem}{Theorem}

\theoremstyle{definition}
\newtheorem*{definition}{Definition}

\newtheorem*{theorem*}{Theorem}

\theoremstyle{remark}
\newtheorem{remark}[equation]{Remark}

\begin{document}

\title[Small Points on Elliptic Curves]
{Global Discrepancy and Small Points on Elliptic Curves}

\date{July 2005}

\author{Matthew Baker}
\email{mbaker@math.gatech.edu}
\address{School of Mathematics, 
         Georgia Institute of Technology, Atlanta, GA 30332-0160}

\author{Clayton Petsche}
\email{clayton@math.uga.edu}
\address{Department of Mathematics, 
         The University of Georgia, Athens, GA 30602-7403}

\begin{abstract}
Let $E$ be an elliptic curve defined over a number field $k$.  
In this paper, we define the ``global discrepancy'' of a finite set 
$Z \subset E(\kbar)$ which in a precise sense measures how far 
the set is from being adelically equidistributed.  We prove an upper bound
for the global discrepancy of $Z$ in terms of the average canonical height of points in $Z$.
We deduce from this inequality a number of consequences.  
For example, we give a new and simple proof of the Szpiro-Ullmo-Zhang
equidistribution theorem for elliptic curves.  We also prove a non-archimedean version of the
Szpiro-Ullmo-Zhang theorem which takes place on the Berkovich analytic space associated to $E$.
We then prove some quantitative `non-equidistribution' theorems for totally real or totally $p$-adic small points.  The
results for totally real points imply similar bounds for points defined over the
maximal cyclotomic extension of a totally real field.  
\end{abstract}

\maketitle


\begin{section}{Introduction}
\label{Introduction}

Let $k$ be a number field, and let $E/k$ be an elliptic curve.
We denote by $\h : E(\kbar) \to \RR$ the N{\'e}ron-Tate canonical height function on $E$, which
is nonnegative and vanishes precisely on the torsion subgroup $E(\kbar)_{\tor}$.
If $Z \subset E(\kbar)$ is a finite set of points, we define the canonical height of $Z$ to be
the average of the canonical heights of the points in $Z$, i.e., 
$\h(Z) = \frac{1}{|Z|} \sum_{P \in Z} \h(P)$.

In $\S$\ref{GlobalDiscrepancySection}, we will define another quantity $\D(Z)$ called the {\em global discrepancy} of $Z$;
it is a positive real number which is defined as a sum of local discrepancies.  
The global discrepancy measures, in a certain precise sense, how far the points of $Z$ are from
being adelically equidistributed.

We will give an upper bound for the global discrepancy $\D(Z)$ in terms of the
canonical height $\h(Z)$, and we will deduce from this inequality various results. 
For example, we will prove an adelic equidistribution theorem whose statement combines ingredients from
both \cite{SzpiroUllmoZhang} and \cite{ChambertLoir}.  In order to state the result, we 
recall that for each place $v \in M_k$, the theory of Berkovich \cite{Berkovich} furnishes an analytic space $E_{\Berk,v}$ which
is compact, Hausdorff, and path-connected, and which contains $E(\CC_v)$ as a dense subspace.  When $v$ is archimedean, 
$E_{\Berk,v} = E(\CC)$, but in the non-archimedean case, $E_{\Berk,v}$ is much larger than $E(\CC_v)$.  For each place
$v$, there is a canonical probability measure $\mu_v$ on $E_{\Berk,v}$ which will be defined in $\S$\ref{BerkovichSection1}. 
With these definitions in mind, we have:

\begin{theorem}
\label{AdelicSUZTheorem}
Let $k$ be a number field, and let $E/k$ be an elliptic curve.  
Fix a place $v \in M_k$ and an embedding $E(\kbar) \hookrightarrow E(\CC_v) \subseteq E_{\Berk,v}$.
Suppose that $\{ P_n \}$ is a sequence of distinct points in $E(\kbar)$ such that $\hhat(P_n) \to 0$,
and let $\delta_n$ be the Borel probability measure on $E_{\Berk,v}$ supported equally on the set $Z_n$ of $\Gal(\kbar/k)$-conjugates of $P_n$.
Then $\delta_n \to \mu_v$ weakly on $E_{\Berk,v}$.
\end{theorem}

The archimedean case of Theorem~\ref{AdelicSUZTheorem} will be established in \S\ref{ArchLocalEquidistributionSection},
and the non-archimedean case in \S\ref{BerkovichSection2}.
The proof of Theorem~\ref{AdelicSUZTheorem} in the archimedean case is more elementary than the 
one given in \cite{SzpiroUllmoZhang}, and is better suited to quantitative refinements.  In the non-archimedean case, the result
generalizes a theorem of \cite{ChambertLoir}, and again our proof is more elementary.
On the other hand, the results of \cite{SzpiroUllmoZhang} and some of the results of \cite{ChambertLoir} 
apply to abelian varieties in general.

\medskip

We will also give a number of quantitative ``non-equidistribution'' results.  Specifically, we will provide 
explicit upper bounds for $|E(L)_{\tor}|$ and lower bounds for 
$\liminf_{P \in E(L)} \h(P)$ and $\inf_{\h(P) \neq 0} \h(P)$
when $L$ is one of the following types of algebraic extensions of $k$:
\begin{itemize}
\item[(i)] $L$ is totally real, i.e., every embedding of $L$ into $\CC$ has image contained in $\RR$.
\item[(ii)] $L$ is totally $p$-adic for some prime number $p$, i.e.,
every embedding of $L$ into $\Qpbar$ has image contained in $\QQ_p$.  
\item[(iii)] $L$ is the maximal cyclotomic extension of $k$ (when $k$ is totally real).
\end{itemize}
The basic idea behind all of these applications of the height-discrepancy inequality is that 
the discrepancy of a finite subset $Z \subset E(L)$ cannot be too small.
It has been previously established (see \cite{Zhang},\cite{Ribet}) that $|E(L)_{\tor}| < \infty$ and
$\liminf_{P \in E(L)} \h(P)>0$ in cases (i) and (iii), so in these cases the novelty in our results 
is that our method of proof leads to completely explicit bounds.

\medskip

We now define the global discrepancy $\D(Z)$ of a finite set $Z = \{ P_1, \ldots P_N \} \subset E(K)$,
where $K/k$ is a finite extension.
The global discrepancy can be thought of as a ``smoothing out'' of the quantity
\[
\Lambda(Z) = \frac{1}{N^2}\sum_{1\leq i,j \leq N} \h(P_i - P_j).
\]
Note that $\Lambda(Z) \geq 0$, and that the parallelogram law furnishes the inequality
\begin{equation}
\label{parlawineq}
\Lambda(Z) \leq 4\h(Z).
\end{equation}
Furthermore, we can decompose $\Lambda(Z)$ as a sum 
\[
\Lambda(Z) = \sum_{v \in M_K} \frac{[K_v : \QQ_v]}{[K:\QQ]} \Lambda_v(Z)
\]
with
\begin{equation}
\label{LocalLambdaEquation}
\Lambda_v(Z) = \frac{1}{N^2} \sum_{\stackrel{1\leq i,j \leq N}{i \neq j}} \lambda_v(P_i - P_j),
\end{equation}
where $\lambda_v$ is an appropriately normalized N{\'e}ron local height function $\lambda_v : E(\CC_v) \backslash \{ O \} \to \RR$.
However, the quantity $\Lambda_v(Z)$ can be negative at archimedean places and at non-archimedean places of bad reduction; this is
closely related to the fact that the singularity of $\lambda_v$ at the origin has forced us to remove the diagonal 
from (\ref{LocalLambdaEquation}).
For the applications we have in mind, it is useful to work with a nonnegative variant of $\Lambda_v(Z)$, and thus of
$\Lambda(Z)$ as well, which can still be bounded explicitly in terms of $\h(Z)$.  At an archimedean place $v$, 
Elkies (see \cite{LangAT}, $\S$VI) accomplishes this via convolution with the heat kernel (see also Faltings \cite{Faltings}).
This gives a one-parameter family $\{ \lambda_t \}_{t > 0}$ of smooth functions $\lambda_t : E(\CC) \to \RR$, defined on all of
$E(\CC)$, such that $\lim_{t \to 0} \lambda_t = \lambda_v$ and 
\[
\frac{1}{N^2} \sum_{1\leq i,j \leq N} \lambda_t(P_i - P_j) > 0.
\]
Although there is no canonical choice for the parameter $t$, the simplest choice 
(which gives nearly optimal estimates) is $t = 1/N$.  We thus define the archimedean local discrepancy 
of a subset $Z \subset E(\CC)$ to be
\[
D_v(Z) = \frac{1}{N^2} \sum_{1\leq i,j \leq N} \lambda_{\frac{1}{N}}(P_i - P_j).
\]

We will see that if $ \{ Z_n \}_{n \geq 1}$ is a sequence of finite subsets of $E(\CC)$ such that
$\lim_{n \to \infty} D_v(Z_n) = 0$, then the sequence of discrete probability measures
$\delta_n$ supported equally on the elements of $Z_n$ converges weakly to the normalized Haar measure $\mu$ on $E(\CC)$.  
The archimedean local discrepancy can thus be thought of as a measure of how far a set $Z$ is from 
being equidistributed in $E(\CC)$.

\medskip

For non-archimedean places $v$, there is a simpler way to modify $\Lambda_v$ in order to obtain a suitable
nonnegative quantity.  As we explain further in $\S$\ref{NonarchLocalDiscrepSection} (see also \cite{ChinburgRumely}), 
for all $P,Q \in E(\CC_v)$ with $P \neq Q$, there is a decomposition
\[
\lambda_v(P - Q) = i_v(P,Q) + j_v(P,Q),
\]
where $i_v(P,Q)$ is a nonnegative arithmetic intersection term which tends to $+\infty$ as $P \to Q$ 
and $j_v(P,Q)$ is a bounded term which has a natural interpretation in terms of the ``skeleton'' 
of $E(\CC_v)$ (see \S\ref{RetractionHomomorphismSection}).
For all $P,Q \in E(\CC_v)$, we set
\[
i_v^*(P,Q) = \left\{ 
\begin{array}{ll}
i_v(P,Q) & P \neq Q \\
0 & P = Q \\
\end{array}
\right.
\]
and 
\[
\lambda_v^*(P-Q) = i_v^*(P,Q) + j_v(P,Q),
\]
and then define the non-archimedean local discrepancy of a set $Z \subset E(\CC_v)$ to be
\[
D_v(Z) = \frac{1}{N^2} \sum_{1\leq i,j \leq N} \lambda_v^*(P_i - P_j).
\]
Note that if $v$ is a place of good reduction, then $D_v(Z) = \Lambda_v(Z)$.

To help justify the definition of the non-archimedean local discrepancy, we will show that for each non-archimedean place
$v \in M_k$, there is a natural probability measure $\mu_v$ on the Berkovich analytic space $E_{\Berk,v}$ such that
if $\{ Z_n \}_{n \geq 1}$ is a sequence of finite subsets of $E(\CC_v)$ with
$\lim_{n \to \infty} D_v(Z_n) = 0$, then the sequence of discrete probability measures
$\delta_n$ supported equally on the elements of $Z_n$ converges weakly to $\mu_v$ on $E_{\Berk,v}$.  
Thus the non-archimedean local discrepancy can also be thought of as a quantitative measure of equidistribution 
(or non-equidistribution) in $E_{\Berk,v}$.

\medskip

Finally, we state the main inequality linking the global discrepancy
\[
\D(Z) = \sum_{v \in M_K} \frac{[K_v : \QQ_v]}{[K:\QQ]} D_v(Z)
\]
of a subset $Z \subset E(\kbar)$ with the canonical height of $Z$.
(Here $K$ is any finite extension of $k$ such that $Z \subset E(K)$ and $E/K$ is semistable; one can show
that $\D(Z)$ does not depend on the choice of $K$.)  Recall that 

\medskip

{\bf Height-Discrepancy Inequality.}
Let $Z=\{P_1, \dots,P_N\}\subseteq E(\kbar)$ be a set of $N$ distinct algebraic points.  Then
\begin{equation}
\label{discrepineq}
\D(Z) \leq 4\h(Z) + \frac{1}{N}\Big(\frac{1}{2}\log N + \frac{1}{12}h(j_E) + \frac{16}{5}\Big),
\end{equation}
where $h(j_E)$ is the logarithmic absolute Weil height of the $j$-invariant $j_E$ of $E/k$.

\medskip

Note that as $N\to\infty$, (\ref{discrepineq}) gives the same asymptotic estimate as (\ref{parlawineq}).
The proof of (\ref{discrepineq}) will be given in Theorem~\ref{mainthm} of \S\ref{GlobalDiscrepancyUpperBoundSection}.

\medskip

The techniques used in this paper combine ideas from various sources.  The use of the local sums $\Lambda_v(Z)$ in the context of 
the equidistribution of small points originates in Baker-Rumely \cite{BakerRumely}.  
A related approach for the multiplicative group $\G_m$ occurs in Bombieri's paper \cite{Bombieri}, which
inspired the Fourier-theoretic discrepancy methods employed here.  Bombieri's approach is a simplification of Bilu's in \cite{Bilu}.
We also note that a series of papers by Hindry-Silverman 
\cite{HindrySilvermanI}, \cite{HindrySilvermanII}, \cite{HindrySilvermanIII}
makes use of explicit lower bounds for the sums $\Lambda_v(Z)$ in order to obtain quantitative results concerning the heights of points in $E(k)$.

\end{section}


\begin{section}{The Archimedean Local Discrepancy}

Throughout this section $E/\C$ is an elliptic curve with j-invariant $j_E$, and $\mu$ denotes the unit Haar measure on the compact group $E(\C)$.  

\begin{subsection}{The Fourier Transform and Laplacian}

Let $\Gamma_E$ denote the dual group of $E(\C)$, 
that is the continuous homomorphisms of $E(\C)$ into the circle group $\T=\{z\in\C\,\mid\,|z|=1\}$.  Given $f\in L^1(E(\C),\mu)$, 
the Fourier transform $\widehat{f}:\Gamma_E\to\C$ is defined by  
\begin{equation*}
\widehat{f}(\gamma)=\int_{E(\C)}f(P)\overline{\gamma(P)}d\mu(P);
\end{equation*}
and similarly, given a signed Borel measure $m$ on $E(\C)$, we define its Fourier-Stieltjes transform $\widehat{m}:\Gamma_E\to\C$ by
\begin{equation*}
\widehat{m}(\gamma)=\int_{E(\C)}\overline{\gamma(P)}dm(P).
\end{equation*}
As $E(\C)$ is compact, its dual $\Gamma_E$ is discrete, and thus each $f\in L^1(E(\C),\mu)$ has a Fourier series
\begin{equation}\label{fs2}
f(P)\sim\sum_{\gamma\in\Gamma_E}\widehat{f}(\gamma)\gamma(P).
\end{equation}
We have equality in $(\ref{fs2})$ provided the right-hand-side is absolutely convergent. 

Denote by $\Sm(E(\C))$ the space of smooth functions on the curve, and by $\Sm'(E(\C))$ its dual, the space of distributions.  
We define the Laplacian on $E(\C)$ as an operator $\Delta:\Sm(E(\C))\to\Sm(E(\C))$ as follows.  Fix a complex uniformization 
\begin{equation}\label{para}
E(\C)\simeq\C/L
\end{equation}
for a normalized lattice $L=\Z+\tau\Z$ ($\tau=a+bi\in\H$), and let $z=x+yi\in\C$ be a complex variable.  
Now, given $g\in\Sm(E(\C))$, we define $\Delta g\in\Sm(E(\C))$ by 
\begin{equation}\label{laplacian}
\Delta g =\frac{b}{2\pi}\Big(\frac{\partial^2}{\partial x^2}+\frac{\partial^2}{\partial y^2}\Big)g.
\end{equation}
It is straightforward to show that this definition does not depend on the choice of $L$ in its homothety class, and thus $\Delta$ is well defined on the curve $E(\C)$.  

Given a uniformization $(\ref{para})$, one can write down explicitly the characters $\gamma\in\Gamma_E$.  
In this case the dual group $\Gamma_E$ is parametrized by the lattice $L$ itself; that is, 
for each lattice point $\omega=n_1+n_2\tau\in L$, we have a character $\gamma_\omega:\C/L\to\T$ given by
\begin{equation}\label{chardef}
\gamma_\omega(z)  = e^{2\pi i (n_1r_1+n_2r_2)},
\end{equation} 
where $r_1,r_2$ are real variables, and $z=x+yi=r_1+r_2\tau\in\C$.  In order to describe the effect of the Laplacian on characters, 
we first define a permutation $\omega\mapsto\omega'$ of the lattice $L$ via the formula
\begin{equation}\label{latperm}
(n_1+n_2\tau)' =n_2-n_1\tau.
\end{equation}
Since $r_2=y/b$ and $r_1=x-ay/b$, applying the Laplacian to the function 
\begin{equation*}
\begin{split}
\gamma_\omega(z)        & = e^{2\pi i(n_1r_1+n_2r_2)} \\
        & = e^{2\pi i(n_1x + b^{-1}(n_2-n_1 a)y)},
\end{split}
\end{equation*}
we find that
\begin{equation}\label{lapchar}
\Delta\gamma_\omega(z)=-\frac{2\pi}{b}|\omega'|^2\gamma_\omega(z).
\end{equation}
\end{subsection}

\begin{subsection}{The N\'eron Function}
\label{ArchNeronFunctionSection}
Let $\lambda:E(\C)\setminus\{O\}\to\R$ denote the N\'eron function, as defined and normalized in \cite{SilvermanII}, $\S$VI.1
and \cite{LangAT}, $\S$II.5.  
The most concrete way to define $\lambda$ is in terms of a complex uniformization $E(\C)\simeq\C/L$, where $L=\Z+\tau\Z$ is a lattice and $\tau=a+bi\in\H$.  
Let $z=x+iy=r_1+r_2\tau$ be a complex variable, with $x,y,r_1,r_2 \in \R$, and let $u=e^{2\pi iz}$ and $q=e^{2\pi i\tau}$.  
In this local coordinate system, the N\'eron function $\lambda:(\C/L)\setminus\{0\}\to\R$ is given by the formula
\begin{equation}\label{exlh}
\lambda(z)=-\frac{1}{2}B_2\Big(\frac{\log|u|}{\log|q|}\Big)\log|q| -\log|1-u| -\sum_{n\geq1}\log|(1-q^nu)(1-q^n/u)|,
\end{equation}
where $B_2(T)=T^2-T+1/6$ is the second Bernoulli polynomial (cf. \cite{SilvermanII}, $\S$VI.3).

A more conceptual definition of $\lambda$ can be given in terms of Green's functions.
First, note that since $\lambda$ is continuous except for a logarithmic singularity at $O$, it is $\mu$-integrable.  
Define distributions $\delta_O, D_\mu, D_{\Delta\lambda} \in \Sm'(E(\C))$ by $\delta_O(g)=g(O)$,
\begin{equation*}
D_{\Delta\lambda}(g)=\int_{E(\C)}\lambda(P)\Delta g(P)d\mu(P),
\end{equation*}
and
\begin{equation*}
D_\mu(g)=\int_{E(\C)}g(P)d\mu(P),
\end{equation*}
where $g\in  \Sm(E(\C))$.

The following proposition, which is proved in \cite{LangAT}, $\S$II, Theorem 5.1,
says that $\lambda$ is the Green's function
on $E(\C)$ with respect to the divisor $(O)$.\footnote{There are different conventions in the literature for what one
means by a Green's function; in Lang's book \cite{LangAT}, for example, the Green's function is twice the N{\'e}ron function $\lambda$, whereas
others take the Green's function to be $-\lambda$.  See the remark on page 22 of \cite{LangAT}.}

\begin{prop}
We have
\begin{equation}\label{g1}
\int_{E(\C)}\lambda(P)d\mu(P)=0,
\end{equation}
and 
\begin{equation}\label{g2}
D_{\Delta\lambda}=D_\mu-\delta_O.
\end{equation}
\end{prop}

More generally, setting $g(P,Q) = g_Q(P) = \lambda(P-Q)$, we have $\int g(P,Q) \, d\mu(P) d\mu(Q)=0$ and
$D_{\Delta g_Q}=D_\mu-\delta_Q$, which means that 
$g(P,Q)$ is the unique normalized Arakelov-Green's function on $E(\CC)$ with
respect to $\mu$ (see \cite{LangAT}, $\S$II, \cite{ChinburgRumely}, and \cite{BakerRumely}).

\medskip

It is evident from $(\ref{lapchar})$ that the characters of $E(\C)$ are eigenfunctions of the Laplacian.  The following proposition shows that the eigenvalues are closely related to the Fourier coefficients of the N\'eron function.  We let $\gamma_0\in\Gamma_E$ denote the trivial character of $E(\C)$.  

\begin{prop}\label{eigenfunction}
Let $\gamma\in \Gamma_E$ be a character of $E(\C)$.  Then
\begin{equation*}
\Delta\gamma = -c_\gamma\gamma,
\end{equation*}
where
\begin{equation*}
c_\gamma = 
\begin{cases}
0 & \text{ if } \gamma=\gamma_0, \\
1/\widehat{\lambda}(\gamma) & \text{ if }  \gamma\neq\gamma_0.
\end{cases} \\
\end{equation*}
Moreover, we have $\widehat{\lambda}(\gamma)\geq0$, with equality if and only if $\gamma=\gamma_0$.
\end{prop}
\begin{proof}
The fact that the characters are eigenfunctions for the Laplacian follows from $(\ref{lapchar})$, from which it also follows 
that $c_\gamma\geq0$, with equality if and only if $\gamma=\gamma_0$.  Applying the distributional identity $(\ref{g2})$ to a nontrivial character $\gamma$, we find $\int_{E(\C)}\lambda(P)\Delta\gamma(P)d\mu(P)=-1$.  
On the other hand,
\begin{equation*}
\begin{split}
\int_{E(\C)}\lambda(P)\Delta \gamma(P)d\mu(P) & = -c_{\gamma}\int_{E(\C)}\lambda(P)\gamma(P)d\mu(P) \\
        & = -c_{\gamma}\widehat{\lambda}(\overline{\gamma}),
\end{split}
\end{equation*}
where $\overline{\gamma}(P)=\overline{\gamma(P)}$ is the complex conjugate of $\gamma\in\Gamma_E$.  Therefore $c_\gamma=1/\widehat{\lambda}(\overline{\gamma})$.  
The proof of the proposition is now complete upon noting that $\widehat{\lambda}(\gamma)=\widehat{\lambda}(\overline{\gamma})$, 
since $\lambda$ is real-valued.
\end{proof}

Given a complex uniformization $(\ref{para})$, we can now give an explicit formula for the Fourier coefficients of the N\'eron function.  
If $\omega\in L$ is a lattice point, and $\gamma_\omega$ is the character given in $(\ref{chardef})$, 
then it follows from Proposition \ref{eigenfunction} and $(\ref{lapchar})$ that
\begin{equation}\label{explicitcoef}
\widehat{\lambda}(\gamma_\omega)=\begin{cases}
0 & \text{ if } \omega=0 \\
\frac{b}{2\pi|\omega'|^2} & \text{ if } \omega\neq0.
\end{cases} 
\end{equation}
\end{subsection}

\begin{subsection}{The Heat Kernel}
The Fourier series 
\begin{equation*}
\lambda(P) \sim \sum_{\gamma\in\Gamma_E}\widehat{\lambda}(\gamma)\gamma(P)
\end{equation*}
for the N\'eron function does not converge absolutely, due to the fact that it has a singularity at $O$.  Following Elkies (cf. \cite{LangAT}, $\S$VI), we sidestep this difficulty by convolving with the heat kernel
\begin{equation*}
g_t(P)=\sum_{\gamma\in\Gamma_E}e^{-t/\widehat{\lambda}(\gamma)}  \gamma(P), 
\end{equation*}
where $t>0$ is a parameter, and put 
\begin{equation}\label{lambdat}
\begin{split}
\lambda_t(P) & = \lambda * g_t(P) \\
        & = \int_{E(\C)}\lambda(Q)g_t(P-Q)d\mu(Q) \\
        & = \sum_{\gamma\in\Gamma_E\setminus\{\gamma_0\}}\widehat{\lambda}(\gamma)e^{-t/\widehat{\lambda}(\gamma)} \gamma(P).
\end{split}
\end{equation}

\end{subsection}

\begin{subsection}{Discrepancy}
Let $Z=\{P_1,\dots,P_N\}\subset E(\C)$ be a set of $N$ points on the elliptic curve $E(\C)$.   We define the discrepancy of this set to be
\begin{equation}\label{adisc}
D(Z)=\frac{1}{N^2}\sum_{1\leq i,j\leq N}\lambda_{\frac{1}{N}}(P_i-P_j).
\end{equation}
The following proposition shows that $D(Z)$ is essentially a classical $L^2$-type discrepancy, measuring the $\mu$-uniform distribution of the set $Z$.

\begin{prop}\label{parsprop}
For each $n\geq1$, let $Z_n\subset E(\C)$ be a set of $N_n$ distinct points, and let $\delta_n$ denote the probability measure on $E(\C)$ 
that assigns a mass of $1/N_n$ at each point of $Z_n$.  If $D(Z_n)\to0$, then the sequence of measures $\delta_n$ converges weakly to Haar measure $\mu$.
\end{prop}
\begin{proof}
In order to show that $\delta_n\to\mu$ weakly, it suffices by Fourier inversion to show that $\widehat{\delta_n}(\gamma)\to\widehat{\mu}(\gamma)=0$ 
for all nontrivial characters $\gamma\in\Gamma_E$.  By Parseval's formula and the fact that $\widehat{\lambda_t}(\gamma_0)=0$, we have 
\begin{equation}\label{napars}
\begin{split}
D(Z_n) & = \sum_{\gamma\in\Gamma_E}\widehat{\lambda_{\frac{1}{N_n}}}(\gamma)|\widehat{\delta_n}(\gamma)|^2 \\
        & = \sum_{\gamma\in\Gamma_E\setminus\{\gamma_0\}}\widehat{\lambda_{\frac{1}{N_n}}}(\gamma)|\widehat{\delta_n}(\gamma)|^2.
\end{split}
\end{equation}
Now, given a nontrivial character $\gamma\in\Gamma_E$ we have $\widehat{\lambda_t}(\gamma)\nearrow \widehat{\lambda}(\gamma)>0$ as $t\to0^+$, by $(\ref{lambdat})$.  It follows that 
\begin{equation}
\begin{split}
\limsup_{n\to\infty}\{|\widehat{\delta_n}(\gamma)|^2/D(Z_n)\} & \leq \limsup_{n\to\infty}\{1/\widehat{\lambda_{\frac{1}{N_n}}}(\gamma)\} \\
        & = 1/\widehat{\lambda}(\gamma) <+\infty.
\end{split}
\end{equation}
The desired limit $\widehat{\delta_n}(\gamma)\to0$ now follows from the assumption that $D(Z_n)\to0$.
\end{proof}

The following result is a quantitative refinement, for the special case of elliptic curves, 
of a general estimate due to N.~Elkies (see \cite{LangAT}, $\S$VI, Theorem 5.1).
Although similar quantitative estimates have been used frequently in the literature 
(see for example \cite{HindrySilvermanII} and \cite{HindrySilvermanIII}), to our knowledge, 
the details underlying these estimates have never been published.  For this reason, we include a proof of
Proposition~\ref{elkies} as an appendix to this paper.  

\begin{prop}\label{elkies}
Let $E/\C$ be an elliptic curve with j-invariant $j_E$, and let 
\begin{equation*}
Z=\{P_1,\dots,P_N\}\subset E(\C)
\end{equation*}
be a set of $N$ distinct points.  Then $D(Z)>0$, and 
\begin{equation}\label{elkies0}
\sum_{\stackrel{1\leq i,j\leq N}{i\neq j}}\lambda(P_i-P_j)\geq N^2D(Z) - \frac{N\log N}{2} -\frac{N}{12}\log^+|j_E| - \frac{16N}{5}.
\end{equation}
\end{prop}

\end{subsection}
\end{section}


\begin{section}{The Non-archimedean Local Discrepancy}
\label{NonarchLocalDiscrepSection}
In this section, $v$ denotes a non-archimedean place of a number field $k$, 
$k_v$ denotes the completion of $k$ with respect to $v$, 
and $\C_v$ denotes the completion of the algebraic closure $\overline{k}_v$.
We also let $\Ovhat$ denote the ring of integers of $\C_v$ and $\Fvbar$ its residue field.
Finally, we let $E/\C_v$ be an elliptic curve with j-invariant $j_E$.

If $|j_E|_v \leq 1$, then $E$ extends to an abelian scheme over $\Ovhat$, i.e., 
there exists a smooth proper model $\E$ over $\Ovhat$ whose special fiber $\bar{\E}$ is 
an elliptic curve over $\Fvbar$.  In this case we say that $E$ has good reduction.

On the other hand, if $|j_E|_v > 1$, then 
by Tate's non-archimedean uniformization theory
(see \cite{SilvermanII} $\S$V, \cite{Robert}, \cite{Tate}), 
there is an analytic isomorphism $E(\C_v) \simeq \C_v^* / q^\ZZ$, where $q \in \C_v^*$ and $|q|_v < 1$.
In this case, we say that $E$ has (split) multiplicative reduction.

\begin{subsection}{The retraction homomorphism}
\label{RetractionHomomorphismSection}

Suppose first that $E$ has good reduction, and let $\E$ be a model for $E$ over $\Ovhat$.
The special fiber $\bar{\E}$ of $\E$ is an elliptic curve over $\Fvbar$, 
and there is a canonical surjective reduction map $\red : E(\CC_v) \to \bar{\E}(\Fvbar)$ which
is a group homomorphism.  Abusing notation somewhat, we define 
\begin{equation}
\label{GoodRed}
\Ebar = \bar{\E}(\Fvbar)
\end{equation}
to be the image of the reduction homomorphism.

\medskip

Now suppose that $E$ has multiplicative reduction.  
If $P \in E(\C_v)$, we write $u(P)$ for the image of $P$ in $\C_v^* / q^\ZZ$ under the Tate map
$E(\C_v) \simeq \C_v^* / q^\ZZ$.
We define the {\em retraction homomorphism} $r : E(\C_v) \to \RR / \ZZ$ by
\begin{equation}
\label{RetractionMap}
r(P) = \frac{\log |u(P)|_v}{\log |q|_v}.
\end{equation}
We call the circle group $\RR / \ZZ$ the {\em skeleton} of $E$; this terminology comes from
Berkovich's theory of analytic spaces and will be discussed further in $\S$\ref{BerkovichSection1}.

Note that since $|\C_v^*|_v = p^{\QQ}$, the image of the retraction homomorphism
is actually contained in the subgroup $\QQ / \ZZ$ of $\RR / \ZZ$.

%

\end{subsection}

\begin{subsection}{The N\'eron Function}
Let $\lambda_v:E(\C_v)\setminus\{O\}\to\R$ denote the N\'eron local height function, 
as defined and normalized in \cite{SilvermanII}, $\S$VI.1.  In order to describe this function explicitly, we let $E_0(\C_v)$ denote the subgroup of $E(\C_v)$ consisting of those points with ``non-singular reduction''.  
More concretely, if $E$ has good reduction then $E_0(\C_v) = E(\C_v)$, and if $E$ has multiplicative
reduction, then $E_0(\C_v) = {\rm ker}(r)$ is the kernel of the retraction homomorphism.

\medskip

\begin{definition}
The N{\'e}ron function is defined as follows:
\begin{itemize}
\item[(i)]
If $P\in E_0(\C_v)\setminus\{O\}$, then
\begin{equation*}
\begin{split}
\lambda_v(P) & = \log^+|z(P)^{-1}|_v + \frac{1}{12}\log^+|j_E|_v \\
        & = \frac{1}{2}\log^+|x(P)|_v + \frac{1}{12}\log^+|j_E|_v,
\end{split}
\end{equation*}
where $x$ and $y$ are the usual coordinate functions associated to a Weierstrass equation for $\E$ over $\Ovhat$, 
and $z=-x/y$ is the standard local parameter at the origin.  
\item[(ii)] 
If $P\in E(\C_v)\setminus E_0(\C_v)$ is a point with singular reduction, then 
\begin{equation*}
\lambda_v(P) = \frac{1}{2}\B_2(r(P))\log^+|j_E|_v,
\end{equation*}
where
\begin{equation}\label{bern}
\begin{split}
\frac{1}{2}\B_2(t) & = \frac{1}{2}(t-[t])^2 - \frac{1}{2}(t-[t]) + \frac{1}{12} \\
        & = \sum_{k\in\Z\setminus\{0\}}\frac{1}{(2\pi k)^2}e^{2\pi ikt}.
\end{split}
\end{equation}
is one-half the periodic second Bernoulli polynomial.
\end{itemize}
\end{definition}

Note that $\log^+|j_E|_v = 0$ when $E$ has good reduction, and 
$\log^+|j_E|_v = -\log |q| > 0$ when $E$ has multiplicative reduction.

\medskip

Since this definition of $\lambda_v$ may appear rather ad-hoc, we offer the following alternative description.
Following \cite{ChinburgRumely}, there is a decomposition 
\begin{equation*}
\lambda_v(P-Q) = i_v(P,Q) + j_v(P,Q),
\end{equation*}
valid for all pairs $(P,Q)\in E(\C_v)\times E(\C_v)$ with $P\neq Q$, where $i_v$ is a local intersection term, and $j_v$ is a function which factors 
through the retraction map $r:E(\C_v)\to \RR / \ZZ$.   

When $E$ has good reduction, 
$i_v(P,Q)$ is (up to a constant multiple) the scheme-theoretic intersection multiplicity of $P$ and $Q$ on $\E$, which in
concrete terms means that 
\[
i_v(P,Q) = \log^+ |z(P-Q)^{-1}|_v
\]

When $E$ has multiplicative reduction, 
one can define $i_v(P,Q)$ using a Tate uniformization via the formula
\[
i_v(P,Q) = \left\{ 
\begin{array}{ll} 
0 & r(P) \neq r(Q) \\
-\log |1 - \frac{u(P)}{u(Q)}|_v & r(P) = r(Q), \\
\end{array}
\right.
\]
where $u(P), u(Q) \in \C_v^*$ are representatives of the classes of $P$ and $Q$, respectively,
in $\C_v^* / q^\ZZ$, chosen in such a way that $|u(P)|_v = |u(Q)|_v$ (which is possible precisely
when $r(P)=r(Q)$).

\medskip

It is easy to see that in both cases, $i_v(P,Q)$ is symmetric in $P$ and $Q$, and
provides a natural measure of how $v$-adically close $P$ and $Q$ are.

\medskip

When $E$ has good reduction, $j_v(P,Q)$ is identically zero.
When $E$ has multiplicative reduction, we define
\[
j_v(P,Q) = \frac{1}{2}\B_2(r(P-Q)) \log^+|j_E|_v.
\]

Note that $j_v(P,Q)$ is also symmetric in $P$ and $Q$, and that
$j_v(P,Q)= \frac{1}{12} \log^+|j_E|_v$ if $r(P)=r(Q)$.  
For those familiar with the language of metrized graphs (see \cite{ChinburgRumely}, \cite{HAOMG}, and \cite{BakerRumely}),
we point out that $g(x,y) = \frac{1}{2}\B_2(r(P-Q))$ is the unique normalized Arakelov-Green's function for
the Haar measure $\mu$ on the circle $\RR / \ZZ$ (and so the function $j_v(P,Q)$ is not as arbitrary as it might appear).

\end{subsection}

\begin{subsection}{The reduction homomorphism}
\label{ReductionHomomorphismSection}

We can use the $v$-adic intersection function $i_v(P,Q)$ to define the reduction $\Ebar$ of $E$ and 
a reduction homomorphism $\red : E \to \Ebar$ without explicit reference to models in the multiplicative
reduction case.

Recall that $i_v(P,Q)$ is a symmetric, nonnegative, real-valued function of $P$ and $Q$, defined for all
$P,Q \in E(\C_v)$ with $P \neq Q$, which measures the $v$-adic proximity of $P$ and $Q$.  
We extend $i_v(P,Q)$ to all $P,Q \in E(\C_v)$ by defining $i_v(P,P)$ to be $+\infty$.

Now define a relation $\sim$ on $E(\C_v)$ by declaring that $P \sim Q$ if and only if
$i_v(P,Q) > 0$.  It is easy to verify that $\sim$ is an equivalence relation.
Intuitively, $P\sim Q$ iff $P$ and $Q$ are congruent modulo the maximal ideal of $\C_v$. 
We therefore {\em define} the reduction $\Ebar$ of $E$ to be the set of equivalence classes
under $\sim$.  If we denote the equivalence class of $P \in E(\C_v)$ by $[P]$, then 
there is an obvious reduction map $\red : E \to \Ebar$ given by $\red(P) = [P]$.
The binary operation given by $[P] + [Q] = [P+Q]$ is easily verified to be well-defined, and
furnishes $\Ebar$ with a natural abelian group structure.  The reduction map from $E$ to $\Ebar$ is then
a group homomorphism.

\medskip

When $E$ has good reduction, and $\E$ is a model over the ring of integers in $\CC_v$ for $E$ whose
special fiber $\bar{\E}$ is smooth, it is easy to see that $\Ebar = \bar{\E}(\Fvbar)$ as abelian groups.

The more interesting case is when $E$ has multiplicative reduction.
In this case, one has the following alternate description of the equivalence relation $\sim$.
By (\ref{RetractionMap}), the kernel of the retraction homomorphism $r : E(\CC_v) \to \QQ/\ZZ$ is 
isomorphic to
\[
\{ u \in \C_v \; : \; |u|=1 \} = \Ovhat^*.
\]
If $\rho : \ker(r) \to \Ovhat^* \to \Fvbar^*$ denotes the isomorphism of $\ker(r)$ with $\Ovhat^*$ followed by the 
natural map from $\Ovhat$ to its residue field, one easily verifies that $P \sim Q$ if and only if $r(P)=r(Q)$ and 
$\rho(P-Q)=1$. 
From this description, it follows easily that there is a commutative diagram of abelian groups
\[
\begin{CD}
\label{EbarExactSequence}
1   @>>>   \Ovhat^*   @>>>   E(\C_v)   @>r>>   \QQ/\ZZ   @>>>   0    \\
@.           @VVV          @VV{\red}V            @|            @.   \\
1   @>>>   \Fvbar^*   @>>>   \Ebar     @>>\bar{r}>   \QQ/\ZZ   @>>>   0    \\
\end{CD}
\]
(Compare with \cite{SGA7}, $\S$11, Expos\'e IX when $\C_v$ is replaced by 
a discretely valued field.)


As {\em sets}, there is a non-canonical isomorphism
\begin{equation}
\label{NonCanonicalIsom} 
\Ebar \simeq \QQ / \ZZ \times \Fvbar^*.  
\end{equation}
Indeed, $\Ebar$ is the disjoint union over all $t \in \QQ / \ZZ$ of 
$\Ebar_t = \{ P \in \Ebar \; : \; \bar{r}(P)=t \}$, and each set $\Ebar_t$ is 
non-canonically isomorphic to $\ker(\bar{r}) = \Fvbar^*$  by translation.
More precisely, each $\Ebar_t$ is a principal homogeneous space for $\Fvbar^*$.

\end{subsection}

\begin{subsection}{Discrepancy}\label{discnonarch}

We can extend the N\'eron function to a function $\lambda_v^*:E(\C_v)\to\R$ defined on the entire curve by setting
\begin{equation}\label{lambdastar}
\lambda_v^*(P) = \begin{cases}
\lambda_v(P) & \text{ if } P\neq O, \\
\frac{1}{12}\log^+|j_E|_v & \text{ if } P=O.
\end{cases} 
\end{equation}
Replacing $\lambda_v$ by $\lambda_v^*$ is in some sense analogous to the convolution $(\ref{lambdat})$ with the heat kernel in the archimedean case.  
One evident difference, however, is that $\lambda_v^*$ is discontinuous, whereas the convolution $(\ref{lambdat})$ smooths out the
singularity at the origin.
If we define $i_v^* : E(\C_v)\times E(\C_v) \to \R$ by the formula
\begin{equation*}
i_v^*(P,Q) =
\begin{cases}
i_v(P,Q) & P \neq Q, \\
0 & P=Q, \\
\end{cases}
\end{equation*}
then 
\begin{equation}\label{decstar}
\lambda_v^*(P-Q) = i_v^*(P,Q) + j_v(P,Q),
\end{equation}
valid for all pairs $(P,Q)\in E(\C_v)\times E(\C_v)$.  


Now, given a set $Z=\{P_1,\dots,P_N\}\subset E(\C_v)$ of $N$ points, we define the local discrepancy $D(Z)$ of this set by
\begin{equation}\label{ndisc}
D(Z) = \frac{1}{N^2}\sum_{1\leq i,j\leq N}\lambda_v^*(P_i-P_j).
\end{equation}
Using $(\ref{decstar})$, we have a decomposition
\begin{equation}\label{discdecomp}
D(Z) =  D_i(Z) + D_j(Z),
\end{equation}
where
\begin{equation}\label{reductiondiscrepancy}
D_i(Z) =  \frac{1}{N^2}\sum_{1\leq i,j\leq N}i_v^*(P_i,P_j)
\end{equation}
is the {\it congruence discrepancy}, and
\begin{equation*}
D_j(Z) =  \frac{1}{N^2}\sum_{1\leq i,j\leq N}j_v(P_i,P_j)
\end{equation*}
is the {\it retraction discrepancy}.

Clearly $D_i(Z) \geq 0$, and we will see in a moment that $D_j(Z) \geq 0$ as well.
Note that if $E$ has good reduction, then $D_j(Z) = 0$ by definition.
So the retraction discrepancy is only relevant in the bad reduction case.

The next two results illustrate the fact that when $E$ has multiplicative reduction, 
the retraction discrepancy $D_j(Z)$ is a measure of the $\mu$-uniform distribution of the retraction $r(Z)$ of the set $Z$, 
where $\mu$ is Haar measure on the circle group $\RR / \ZZ$.

\begin{prop}\label{retprop}
Let $Z=\{P_1,\dots,P_N\}\subset E(\C_v)$ be a set of $N$ distinct points.
If $E$ has multiplicative reduction, we have
\begin{equation}\label{reddisc}
D_j(Z) = {\log^+|j_E|_v}\sum_{k\in\Z\setminus\{0\}} \frac{1}{(2\pi k)^2} \bigg|\sum_{j=1}^N \frac{e^{2\pi ikr(P_j)}}{N}\bigg|^2 \geq 0.
\end{equation}
\end{prop}
\begin{proof}
The formula for $D_j(Z)$ follows immediately from Parseval's formula and the Fourier expansion $(\ref{bern})$ of the second Bernoulli polynomial.  
The fact that $D_j(Z)$ is positive is then an easy consequence of the explicit formula.
\end{proof}

In addition, one has the following analogue for the circle group $\RR / \ZZ$ of Proposition~\ref{parsprop}.
For the statement, we define the discrepancy of a Borel probability measure $\nu$ on $\RR / \ZZ$ to be
\[
D(\nu) = \iint f(x-y) d\nu(x)d\nu(y),
\]
where $f(t) = \frac{\log^+|j_E|}{2}\B_2(t)$.  Letting $\widehat{\nu}(k) = \int_{\RR / \ZZ} e^{-2\pi i k t} d\nu(t)$ and 
\begin{equation*}
\widehat{f}(k) =  \int_{\RR / \ZZ} f(t)e^{-2\pi i k t} dt = 
\begin{cases}
0 & \text{ if } k=0, \\
\frac{\log^+|j_E|_v}{(2\pi k)^2} & \text{ if }  k\neq0,
\end{cases} 
\end{equation*}
we then have
\begin{equation}
\label{circlemeasurediscrep}
D(\nu) = \sum_{k \in \ZZ} \hat{f}(k)|\widehat{\nu}(k)|^2
\end{equation}
by Parseval's formula.

If $E$ has multiplicative reduction and $Z = \{ P_1,\ldots,P_N \} \subset E(\C_v)$ is a set of $N$ distinct points,
and if we define
\[
\delta_Z = \frac{1}{N} \sum_{i=1}^N \delta_{r(P_i)}
\]
to be the natural discrete measure supported on the multiset $r(Z)$,
then it follows from the definitions that $D_j(Z) = D(\delta_Z)$.

\begin{prop}
\label{circleprop}
If $\nu_n$ is a sequence of probability measures converging weakly on $\RR / \ZZ$ to $\nu$, then 
$D(\nu_n) \to D(\nu)$.  Furthermore, $D(\nu) = 0$ if and only if $\nu$ is the normalized
Haar measure $\mu$ on $\RR / \ZZ$.
\end{prop}

\begin{proof}
The convergence of discrepancies follows from the continuity of $f(x-y)$ and the fact that
$\nu_n \times \nu_n$ converges weakly to $\nu \times \nu$ on $\RR / \ZZ \times \RR / \ZZ$.  The last statement
follows from (\ref{circlemeasurediscrep}) and the fact that 
$\hat{f}(k) > 0$ for all $k \in \ZZ\setminus\{0\}$,
since $\mu$ is characterized by the fact that $\widehat{\mu}(k)=0$ for all 
$k \in \ZZ\setminus\{0\}$.
\end{proof}

Finally, we have the following result, which is a non-archimedean analogue of Proposition \ref{elkies}. 
\begin{prop}\label{elkiesna}
Let $Z=\{P_1,\dots,P_N\}\subset E(\C_v)$ be a set of $N$ distinct points.  Then $D(Z)\geq0$ and 
\begin{equation}\label{elkiesnaineq}
\sum_{\stackrel{1\leq i,j\leq N}{i\neq j}}\lambda_v(P_i-P_j) = N^2D(Z) - \frac{N}{12}\log^+|j_E|_v.
\end{equation}
\end{prop}
\begin{proof}
The non-negativity of the discrepancy follows from $(\ref{discdecomp})$ and the fact that $D_i(Z)$ and $D_j(Z)$ are non-negative.  
The identity $(\ref{elkiesnaineq})$ follows at once from the definition $(\ref{lambdastar})$ of $\lambda_v^*$. 
\end{proof}
\end{subsection}
\end{section}


\begin{section}{Global Discrepancy}
\label{GlobalDiscrepancySection}

Let $k$ be a number field, and let $E/k$ be an elliptic curve with j-invariant $j_E$. 

\begin{subsection}{Definition of the global discrepancy}

In this section, we define the global discrepancy $D(Z)$ of a set $Z=\{P_1, \dots,P_N\}\subseteq E(\kbar)$ of algebraic points.  

Let $K$ be any number field such that $Z\subset E(K)$ and $E/K$ is semistable.  If $v$ is a place of $K$, then $E$ is defined over $K_v$, 
and we can view $Z$ as a subset of $E(\C_v)$ via the embedding $K\hookrightarrow\C_v$.  
We let $D_v(Z)$ denote either the complex local discrepancy defined in $(\ref{adisc})$, 
or the non-archimedean local discrepancy defined in $(\ref{ndisc})$.  We now define the global discrepancy 
\begin{equation}\label{gd}
\D(Z)=\sum_{v\in\M_K}\frac{[K_v:\Q_v]}{[K:\Q]}D_v(Z).
\end{equation}
It is straightforward to show that this quantity is an absolute diophantine invariant of the data $k$, $E$, and $Z$, and doesn't depend on the choice of the number field $K$.  

Finally, we recall the definition of the logarithmic absolute Weil height $h:\overline{\Q}\to[0,+\infty)$: given $\alpha\in\overline{\Q}$, select a number field $K$ containing $\alpha$, and define
\begin{equation*}
h(\alpha) = \sum_{v\in M_K}\frac{[K_v:\Q_v]}{[K:\Q]}\log^+|\alpha|_v.
\end{equation*}
\end{subsection}
Again, as is well known, the value of $h(\alpha)$ does not depend on the choice of $K$.

\begin{subsection}{An upper bound for the global discrepancy}
\label{GlobalDiscrepancyUpperBoundSection}
Given a set of algebraic points $Z=\{P_1, \dots,P_N\}\subseteq E(\kbar)$, define its height 
\begin{equation*}
\h(Z)=\frac{1}{N}\sum_{j=1}^{N}\h(P_j)
\end{equation*}
to be the average of the height of its points.

\begin{thm}\label{mainthm}
Let $Z=\{P_1, \dots,P_N\}\subseteq E(\kbar)$ be a set of $N$ distinct algebraic points.  Then
\begin{equation}\label{mainbound}
\D(Z) \leq 4\h(Z) + \frac{1}{N}\Big(\frac{1}{2}\log N + \frac{1}{12}h(j_E) + \frac{16}{5}\Big).
\end{equation}
\end{thm}
\begin{proof}
Let $K$ be any number field such that $Z\subset E(K)$ and $E/K$ is semistable.  Let $d=[K:\Q]$, and for each place $v$ of $K$, 
let $d_v=[K_v:\Q_v]$ and
\begin{equation}\label{sv}
\Lambda_v(Z)=\frac{1}{N^2}\sum_{\stackrel{1\leq i,j\leq N}{i\neq j}}\lambda_v(P_i-P_j).
\end{equation}
If $v$ is an archimedean place, then
\begin{equation*}
\Lambda_v(Z)\geq D_v(Z) - \frac{\log N}{2N} -\frac{1}{12N}\log^+|j_E|_v - \frac{16}{5N}
\end{equation*}
by Proposition \ref{elkies}.  If $v$ is a non-archimedean place, then 
\begin{equation}\label{sigmav}
\Lambda_v(Z) = D_v(Z) -\frac{1}{12N}\log^+|j_E|_v
\end{equation}
by Proposition \ref{elkiesna}.  Summing over all places of $K$, we have
\begin{equation*}
\begin{split}
\Lambda(Z) & = \sum_{v\in M_K}\frac{d_v}{d}\Lambda_v(Z) \\
        & \geq \D(Z) - \frac{\log N}{2N} -\frac{1}{12N}\sum_{v\in M_K}\frac{d_v}{d}\log^+|j_E|_v - \frac{16}{5N} \\
        & = \D(Z) - \frac{\log N}{2N} -  \frac{1}{12N}h(j_E)-\frac{16}{5N},
\end{split}
\end{equation*}
where we have used the fact that $\sum_{v\mid\infty}d_v=d$.

On the other hand, by the parallelogram law, and the fact that the N\'eron-Tate height is nonnegative, we have $\h(P_i-P_j)\leq 2\h(P_i) +2\h(P_j)$ 
for all $1\leq i,j\leq N$.   In view of this and the decomposition
\begin{equation}\label{decomp}
\h(P)=\sum_{v\in M_K}\frac{d_v}{d}\lambda_v(P),
\end{equation}
valid for all $P\in E(K)\setminus\{O\}$ (cf. \cite{SilvermanII}, $\S$VI.2), we have the upper bound 
\begin{equation}\label{globalub}
\begin{split}
\Lambda(Z) & = \sum_{v\in M_K}\frac{d_v}{d}\Lambda_v(Z) \\
        & = \frac{1}{N^2}\sum_{\stackrel{1\leq i,j\leq N}{i\neq j}}\h(P_i-P_j) \\
        & \leq \frac{1}{N^2}\sum_{\stackrel{1\leq i,j\leq N}{i\neq j}}(2\h(P_i) +2\h(P_j)) \\
        & \leq \frac{4}{N}\sum_{1\leq i\leq N}\h(P_i) \\
        & = 4\h(Z).
\end{split}
\end{equation}
Combining the upper and lower bounds for the sum $\Lambda(Z)$, we deduce $(\ref{mainbound})$.
\end{proof}
\end{subsection}
\end{section}

In $\S$\ref{LocalEquidistributionSection}, we will use Theorem~\ref{mainthm} in conjunction with
the following inequality relating local and global discrepancies:

\begin{prop}
\label{LocalGlobalInequalityProp}
Let $v_0 \in M_k$ be a place of $k$, and let $Z$ be a finite subset of $E(\kbar)$.
Then $D_{v_0}(Z) \leq [k:\QQ] \D(Z)$.
\end{prop}

\begin{proof}
Select a Galois extension $K/k$ containing $k(Z)$, and over which $E$ has everywhere semistable reduction.
Since $Z$ is $\Gal(K/k)$-stable, we have $D_{v_0}(Z)=D_{w}(Z)$ for all places $w\in\M_K$ that are $\Gal(K/k)$-conjugates of $v_0$.  Thus
\begin{equation}
\begin{split}
\D(Z) & = \sum_{v\in\M_K}\frac{[K_v:\Q_v]}{[K:\Q]}D_v(Z) \\
        & \geq \sum_{\sigma \in\Gal(K/k)}\frac{[K_{\sigma v_0}:\Q_{\sigma v_0}]}{[K:\Q]}D_{\sigma v_0}(Z) \\
        & = \Big\{\sum_{\sigma \in\Gal(K/k)}\frac{[K_{\sigma v_0}:\Q_{\sigma v_0}]}{[K:\Q]}\Big\}D(Z) \\
        & = [k:\Q]^{-1}D(Z).
\end{split}
\end{equation}
\end{proof}


\begin{section}{The Local Equidistribution of Small Points}
\label{LocalEquidistributionSection}

\begin{subsection}{The Szpiro-Ullmo-Zhang equidistribution theorem for elliptic curves}
\label{ArchLocalEquidistributionSection}

In this section, we show how one can obtain from Theorem~\ref{mainthm} an ``elementary'' (in the sense that it does not use
Arakelov intersection theory) proof of the Szpiro-Ullmo-Zhang equidistribution theorem \cite{SzpiroUllmoZhang} for elliptic curves.

\begin{thm}[Szpiro-Ullmo-Zhang]
Let $k$ be a number field, and let $E/k$ be an elliptic curve.  Fix an embedding $\kbar\hookrightarrow\CC$, so 
that $E(\kbar) \subset E(\CC)$.
Suppose that $\{ P_n \}$ is a sequence of distinct points in $E(\kbar)$ such that $\hhat(P_n) \to 0$.
Let $\delta_n$ be the Borel probability measure on $E(\CC)$ supported equally on the set $Z_n$ of $\Gal(\kbar/k)$-conjugates of $P_n$.
Then $\delta_n \to \mu$ weakly on $E(\CC)$, where $\mu$ denotes the normalized Haar measure on $E(\CC)$.
\end{thm}

\begin{proof}
By Proposition~\ref{parsprop}, it suffices to show that $D(Z_n)\to0$, where $D(Z_n)$ denotes the archimedean local 
discrepancy of $Z_n\subset E(\C)$.  
Let $\D(Z_n)$ denote the global discrepancy of the set $Z_n\subset E(\kbar)$.  Then $D(Z_n)\leq[k:\Q]\D(Z_n)$ by 
Proposition~\ref{LocalGlobalInequalityProp}.
Since the points $P_n$ are distinct and $\h(P_n)\to0$, it follows from the Northcott finiteness principle 
that $|Z_n|\to+\infty$.  (The Northcott finiteness principle is the statement 
that for any $M>0$, there are only finitely many points $P$ with given degree and $\h(P) \leq M$.) 
Also, $\h(Z_n) = \h(P_n) \to 0$ by the Galois-invariance of the canonical height.
Therefore, from Theorem~\ref{mainthm} we deduce that $D(Z_n) \leq [k:\Q]\D(Z_n) \to0$, which completes the proof.
\end{proof}

\end{subsection}

\begin{subsection}{The Berkovich analytic space associated to an elliptic curve}\label{BerkovichSection1}

In the next section, we will prove a non-archimedean version of the Szpiro-Ullmo-Zhang theorem for elliptic curves.  The natural
context for the result is the Berkovich analytic space attached to an elliptic curve, 
so we first recall some basic facts from \cite{Berkovich}, $\S$IV (see also \cite[\S7.2]{FvdP} and \cite{Thuillier}).
To make our results more accessible, we summarize all of the properties of Berkovich analytic spaces needed to understand 
our proof.

Let $E$ be an elliptic curve defined over $\CC_v$.
The Berkovich analytic space $E_{\Berk} = E_{\Berk,v}$ over $\C_v$ associated to $E$ is a compact, Hausdorff, and 
path-connected topological space which
contains $E(\C_v)$ (endowed with its usual ultrametric topology) as a dense subspace.
We define the {\em skeleton} $\Sigma$ of $E_{\Berk}$ (c.f.~$\S$\ref{NonarchLocalDiscrepSection})
to be a single point when $E$ has good reduction, and to be the circle $\RR / \ZZ$ when $E$ has
multiplicative reduction.
The skeleton $\Sigma$ can naturally be viewed as a subspace of $E_{\Berk}$.
Moreover, there is a strong deformation retraction $r : E_{\Berk} \to \Sigma$ which extends the
``retraction homomorphism'' discussed in $\S$\ref{NonarchLocalDiscrepSection}.
In particular, when $E$ has good reduction, the space $E_{\Berk}$ is contractible, and when $E$ has multiplicative reduction,
its fundamental group is isomorphic to $\ZZ$.

Each connected component of $E(\C_v) \backslash \Sigma$ is homeomorphic to an open Berkovich disk, which can be
given the structure of an infinite real tree (or, in Berkovich's terminology,
a simply-connected one-dimensional quasipolyhedron), see \cite{RumelyBaker}, \cite{Berkovich}, \cite{Thuillier}, 
and \cite{FavreJonsson}.
We may thus think of $E_{\Berk}$ as a family of infinite trees glued together along the skeleton $\Sigma$.

\begin{figure}[!ht]
\scalebox{.2}{\includegraphics{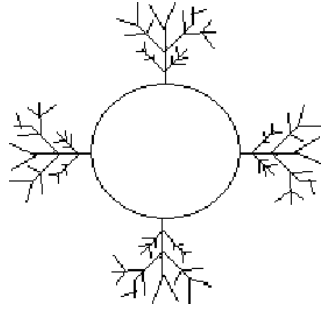}} \caption{The Berkovich space associated to an elliptic curve with multiplicative reduction.}
\label{EllipticBerk}
\end{figure}

More precisely, the complement of $\Sigma$ in $E_{\Berk}$ is a disjoint union of connected open sets
\[
E_{\Berk} \backslash \Sigma = \amalg_{\alpha \in \Ebar} \BB_{\alpha},
\]
with each $\BB_{\alpha}$ homeomorphic to an open Berkovich disk.
We have $\BB_{\alpha} \cap E(\C_v)= B_{\alpha}$, where
$B_{\alpha} = \{ P \in E(\C_v) \; : \;  \red(P) = \alpha \}$.
When $E$ has multiplicative reduction, 
the closure of $\BB_\alpha$ in $E_{\Berk}$ is $\BB_\alpha \cup \{ \bar{r}(\alpha) \}$,
where we identify $\bar{r}(\alpha) \in \QQ / \ZZ$ with its image in $\Sigma \cong \RR / \ZZ \subset E_{\Berk}$.
The set of balls $\BB_{\alpha}$ whose closure meets a given point $t \in \Sigma$ is either empty 
(if $t \not\in \QQ / \ZZ$) or else forms a principal homogeneous space for $\Fvbar^*$ (if $t \in \QQ / \ZZ$).

\medskip

We now define a canonical probability measure $\mu = \mu_v$ (analogous to normalized Haar measure in the archimedean case) 
which is supported on $\Sigma \subset E_{\Berk}$.

\medskip

\begin{definition}
The measure $\mu$ is defined as follows:
\begin{itemize}
\item[(i)] If $E$ has good reduction, $\mu$ is the Dirac measure concentrated on the skeleton $\Sigma$, which consists of a single point.
\item[(ii)] If $E$ has multiplicative reduction, $\mu$ is the uniform probability measure (i.e., Haar measure) 
supported on the circle $\Sigma \cong \RR / \ZZ$.
\end{itemize}
\end{definition}

\begin{remark}
There is a natural way to define a Laplacian operator on $E_{\Berk}$ (see \cite{Thuillier}), and also
to extend the function $g(P,Q) = \lambda_v(P-Q)$ from $E(\C_v)$ to $E_{\Berk}$.
Having done so, one can show that $g(P,Q)$ is the unique normalized Arakelov-Green's function on $E_{\Berk}$ 
with respect to the measure $\mu$ (compare with \S\ref{ArchNeronFunctionSection}).
\end{remark}

For each $\alpha \in \Ebar$, fix a point $Q_\alpha \in E(\CC_v)$ with $\red(Q_\alpha) = \alpha$.
For $0 < r < 1$, let $\BB^+(\alpha,r) \subset \BB_{\alpha}$ be the closure in $E_{\Berk}$ of the set
\[
B^+(\alpha,r) = \{ P \in E(\CC_v) \; : \; i_v(P,Q_\alpha) \geq \log(\frac{1}{r}) \}.
\]

Then 
\[
\BB_{\alpha} = \cup_{r < 1} \BB^+(\alpha,r),
\]
and we have
\begin{equation}
\label{ivinequality}
P,Q \in B^+(\alpha,r) \Rightarrow i_v(P,Q) \geq \log(\frac{1}{r}).
\end{equation}

\end{subsection}

\begin{subsection}{A Berkovich space analogue of the Szpiro-Ullmo-Zhang theorem for elliptic curves}\label{BerkovichSection2}

In this section, we prove our non-archimedean version of the Szpiro-Ullmo-Zhang theorem for elliptic curves.  

\medskip

Let $Z \subset E(\C_v)$ be a subset of cardinality $N = N(Z)$.
For each $\alpha \in \Ebar$ and each $0 < r < 1$, define 
\[
N_{\alpha,r}(Z) = | Z \cap B^+(\alpha,r) |.
\]
Also, define
\[
D_{i,r}(Z) = \sum_{\alpha \in \Ebar} \left( \frac{N_{\alpha,r}(Z)}{N(Z)}  \right)^2.
\]

The following result, which relates $D_{i,r}(Z)$ to the congruence discrepancy $D_i(r)$ defined
in (\ref{reductiondiscrepancy}), is in some ways analogous to Proposition~\ref{elkies}.

\begin{prop}
\label{reductdiscrepprop}
Let $Z=\{P_1,\dots,P_N\}\subset E(\C_v)$ be a set of $N$ distinct points, and let $0<r<1$.  Then
\begin{equation}\label{epdisc}
D_{i,r}(Z) \leq \frac{D_i(Z)}{\log(\frac{1}{r})} + \frac{1}{N}.
\end{equation}
\end{prop}

\begin{proof}
Recall from (\ref{ivinequality}) that if
$P,Q \in B^+(\alpha,r)$, then $i_v(P,Q) \geq \log(\frac{1}{r})$.  
By definition, there are $N_{\alpha,r}(Z)$ points of $Z$ in 
$B^+(\alpha,r)$, and therefore there are
$N_{\alpha,r}(Z) (N_{\alpha,r}(Z) - 1)$ pairs $(P,Q)$ with $P,Q \in Z \cap B^+(\alpha,r)$ and $P \neq Q$.
It follows that
\begin{equation}\label{epdisc1}
\begin{split}
N^2 D_i(Z) &= \sum_{\stackrel{P,Q \in Z}{P \neq Q}}i_v(P,Q) \\
&\geq \log(1/r)\sum_{\alpha \in \Ebar} N_{\alpha,r}(Z) (N_{\alpha,r}(Z) - 1), \\
\end{split}
\end{equation}
and thus
\begin{equation}\label{epdisc2}
\begin{split}
\frac{N^2 D_i(Z)}{\log(\frac{1}{r})} &\geq N^2 D_{i,r}(Z) - \sum_{\alpha \in \Ebar} N_{\alpha,r}(Z) \\
&\geq N^2 D_{i,r}(Z) - N, \\
\end{split}
\end{equation}
which gives the desired result.
\end{proof}

For each $n\geq 1$, let $Z_n \subset E(\C_v)$ be a set consisting of $N_n$ distinct points, and
let $\delta_n$ denote the probability measure on $E_{\Berk}$ supported equally at the elements of $Z_n$.
Note that by definition, we have
\[
\delta_n(\BB^+(\alpha,r)) = \frac{N_{\alpha,r}(Z_n)}{N(Z_n)}.
\]

\begin{prop}
\label{MeasureTheoryProp2}
\begin{itemize}
\item[(a)] Suppose $0 < r < 1$.
If $\lim_{n\to\infty} D_{i,r}(Z_n) = 0$, then 
\[
\lim_{n\to\infty} \frac{ N_{\alpha,r}(Z_n)}{N(Z_n)} = 0
\]
for all $\alpha \in \Ebar$
\item[(b)] Suppose that for all $\alpha \in \Ebar$ and all $0 < r < 1$, we have
\[
\lim_{n\to\infty} \frac{N_{\alpha,r}(Z_n)}{N(Z_n)} = 0.
\]
Then every subsequential limit $\nu$ of $\{ \delta_n \}$ is supported on the skeleton $\Sigma \subset E_{\Berk}$.
\end{itemize}
\end{prop}

\begin{proof}

Part (a) follows from the fact that
\[
\begin{aligned}
0 &\leq \liminf_{n\to\infty} \left( \frac{N_{\alpha,r}(Z_n)}{N(Z_n)} \right)^2 \\
&\leq \limsup_{n\to\infty} \left( \frac{N_{\alpha,r}(Z_n)}{N(Z_n)} \right)^2 \\
&\leq \limsup_{n\to\infty} D_{i,r}(Z_n) = 0. \\
\end{aligned}
\]

To prove (b), fix $\alpha \in \Ebar$ and $0 < r < 1$, and let $A  = \BB^+(\alpha,r)$. 
Let $\nu$ be any weak subsequential limit of the sequence $\{ \delta_n \}$.

We claim that $\nu(A) = 0$.  
To see this, fix $r'$ with $r < r' < 1$ and let $B = \BB^+(\alpha,r')$.
Since the closure of the complement of $B$ is disjoint from $A$,
it follows from Urysohn's Lemma (which is valid for every compact Hausdorff topological space)
that we can choose a continuous function $f$ on $E_{\Berk}$ such that $f \equiv 1$ on $A$, $f \equiv 0$ outside $B$, 
and $0 \leq f \leq 1$ on $B$.  Then
\begin{equation}
\label{eqn1}
\delta_n(A)  \leq \int_{E_{\Berk}} f \, d\delta_n \leq  \delta_n(B) 
\end{equation}
and
\begin{equation}
\label{eqn2}
\nu(A)  \leq  \int_{E_{\Berk}} f \, d\nu  \leq  \nu(B). 
\end{equation}

By hypothesis, $\lim_{n \to \infty} \delta_n(A) = \lim_{n \to \infty} \delta_n(B) = 0$.  
Therefore 
\[
\lim_{n \to \infty} \int_{E_{\Berk}} f \, d\delta_n = 0
\]
by (\ref{eqn1}).
Since $\lim_{n \to \infty} \int_{E_\Berk} f \, d\delta_n = 
\int_{E_{\Berk}} f \, d\nu$ by weak convergence, we have 
$\int_{E_{\Berk}} f \, d\nu = 0$ and therefore $\nu(A) = 0$ by (\ref{eqn2}).  
This proves the claim.

Thus $\nu(\BB_{\alpha}) = 0$ for all $\alpha \in \Ebar$, since 
$\BB_{\alpha} = \cup_{r < 1} \BB^+(\alpha,r)$.

Since $\Ebar$ is countable by (\ref{NonCanonicalIsom}), we have 
$\nu(U) = 0$ as well, where $U = \cup_{\alpha \in \Ebar} \BB_\alpha$,
so that $U$ is contained in the complement of the support of $\nu$.
But $\cup_{\alpha \in \Ebar} \BB_\alpha = E_{\Berk} \backslash \Sigma$,
so $\nu$ is supported on $\Sigma$ as claimed.  
\end{proof}

From these results, we deduce the following non-archimedean analogue of Proposition~\ref{parsprop} for
the local discrepancy defined by (\ref{ndisc}):

\begin{cor}
\label{ReductionDiscrepancyCor}
If $\lim_{n \to \infty} D(Z_n) = 0$, then $\delta_n$ converges weakly to $\mu$ on $E_{\Berk}$.
\end{cor}

\begin{proof}
Since $E_{\Berk}$ is compact, after passing to a subsequence if necessary, it suffices to prove that if $\delta_n \to \nu$ then $\nu = \mu$.
By (\ref{discdecomp}), there is a decomposition $D(Z_n) =  D_i(Z_n) + D_j(Z_n)$ into nonnegative terms,
and thus $\lim_{n \to \infty} D_i(Z_n) = D_j(Z_n) = 0$.
Since $D_i(Z_n) \geq 0$ for all $n$, and since $\lim_{n \to \infty} D_i(Z_n) = 0$, Proposition~\ref{reductdiscrepprop} implies that
$\lim_{n\to\infty} D_{i,r}(Z_n) = 0$ for every $0 < r < 1$ and every $\alpha \in \Ebar$. 
It follows from assertions (a) and (b) of Proposition~\ref{MeasureTheoryProp2} that $\nu$ is supported on $\Sigma$.
If $\Sigma$ consists of a point, then we're done.  
In the case where $\Sigma$ is a circle, Proposition~\ref{circleprop} and the fact that $\lim_{n\to\infty} D_j(Z_n)  = 0$ show that $D(\nu) = 0$ 
and hence $\nu = \mu$.
\end{proof}

We apply these observations to obtain the following non-archimedean version of the Szpiro-Ullmo-Zhang equidistribution theorem
for elliptic curves, which generalizes a result of Chambert-Loir \cite{ChambertLoir}) (see also \cite{BakerRumely}).

\begin{thm}
\label{BerkovichSUZTheorem}
Let $k$ be a number field, and let $E/k$ be an elliptic curve.  
Fix a non-archimedean place $v$ of $k$, and an embedding of $\kbar$ into $\CC_v$ (this allows us to consider $E(\kbar)$ as a subset of 
$E_{\Berk,v}$).
Suppose that $\{ P_n \}$ is a sequence of distinct points in $E(\kbar)$ such that $\hhat(P_n) \to 0$.
Let $\delta_n$ be the Borel probability measure on $E_{\Berk,v}$ supported equally on the set $Z_n$ of $\Gal(\kbar/k)$-conjugates of $P_n$.
Then $\delta_n \to \mu_v$ weakly on $E_{\Berk,v}$.
\end{thm}

\begin{proof}
By Theorem~\ref{mainthm} and Proposition~\ref{LocalGlobalInequalityProp}, our hypotheses imply that the local discrepancy
$D(Z_n)$ tends to $0$.  Thus $\delta_n \to \mu_v$ by Corollary~\ref{ReductionDiscrepancyCor}.
\end{proof}

\begin{remark}
For elliptic curves with multiplicative reduction, Chambert-Loir \cite{ChambertLoir} does not formulate his results
in terms of the Berkovich analytic space $E_{\Berk}$; rather, he proves (in the terminology of \S\ref{discnonarch}) 
that $\delta_{r(Z_n)}$ converges weakly to $\mu$ on the circle group $\Sigma$.  In this sense, our result is more general
than his.  On the other hand, in the good reduction case, Chambert-Loir proves a generalization of 
Theorem~\ref{BerkovichSUZTheorem} to the Berkovich space attached to an arbitrary abelian variety.
\end{remark}

\end{subsection}

\end{section}


\begin{section}{Quantitative Applications}

The general idea behind the results in this section is that,
under certain special conditions, a set of global points cannot be too uniformly distributed at all places.  
Specifically, we will use Theorem~\ref{mainthm} to deduce quantitative upper bounds for the number, and lower bounds for the height, 
of such points.  We begin by establishing two useful lemmas.

\begin{lem}\label{tlem}
Let $N\geq1$ satisfy the bound $N\leq A\log N+B$ for constants $A>0,B\geq0$.  
Then $N\leq (e/(e-1))(A\log A+B) < \frac{8}{5}(A\log A + B)$. 
\end{lem}
\begin{proof}
Let $x>1$ and $y>e$ be real numbers, related by the equation $y=x/\log x$.  Note that
\begin{equation*}
\begin{split}
y\log y & = \frac{x}{\log x}(\log x-\log\log x) \\
        & = x\Big(1-\frac{\log\log x}{\log x}\Big) \\
        & \geq x\Big(1-\frac{1}{e}\Big),
\end{split}
\end{equation*}
and therefore $x\leq(\frac{e}{e-1})y\log y$.  We will apply this inequality with $x=Ne^{B/A}>1$.  By our assumption that $N\leq A\log N+B$, 
we have $Ne^{B/A}\leq Ae^{B/A}\log(Ne^{B/A})$, which implies that $y=x/\log x \leq Ae^{B/A}$.  We then have
\begin{equation*}
\begin{split}
Ne^{B/A}        & = x \\
        & \leq \Big(\frac{e}{e-1}\Big)y\log y \\
        & \leq \Big(\frac{e}{e-1}\Big)Ae^{B/A}\log(Ae^{B/A}) \\
        & = \Big(\frac{e}{e-1}\Big)e^{B/A}(A\log A + B),
\end{split}
\end{equation*}
which concludes the proof.
\end{proof}

\begin{lem}\label{smalllem}
Let $L$ be a subfield of $\Qbar$, and let $E/L$ be an elliptic curve.  If 
\begin{equation*}
|\{P\in E(L)\,\mid\, \h(P)\leq S\}| \leq T 
\end{equation*}
for constants $S\geq0$ and $T\geq1$, then 
\begin{equation*}
\begin{split}
(i) & \,\,\, |E(L)_{\tor}| \leq T,  \\
(ii) & \,\,\, \liminf\{\h(P)\,\mid\,P\in E(L)\} \geq S,\text{ and}  \\
(iii) & \,\,\, \h(P)\geq S/T^2\text{ for all non-torsion points }P\in E(L). 
\end{split}
\end{equation*}

\end{lem}
\begin{proof}
The bounds $(i)$ and $(ii)$ are trivial.  In order to see $(iii)$, consider a non-torsion point $P\in E(L)$, and let $M\geq1$ be the largest 
integer such that $\h((M-1)P)\leq S$.  Thus $\h(jP)\leq S$ for all $0\leq j\leq M-1$, and it follows from the hypothesis that $M\leq T$.  
By the maximality of $M$ we have $M^2\h(P)=\h(MP)>S$, and therefore $\h(P) > S/M^2 \geq S/T^2$.
\end{proof}

\begin{subsection}{Small Totally Real Points}\label{app1}

\begin{thm}\label{trtheorem}
Let $\Q_{\tr}$ denote the maximal totally real subfield of $\Qbar$, and let $E/\Q_{\tr}$ be an elliptic curve with j-invariant $j_E$.  Then 
\begin{equation}\label{trtor}
|E(\Q_{\tr})_{\tor}| \leq 3\,h^*(j_E)^2;
\end{equation}
\begin{equation}\label{trliminf}
\liminf\{\h(P)\,\mid\, E(\Q_{\tr})\} \geq \frac{1}{24\,h^*(j_E)};
\end{equation}
and if $P\in  E(\Q_{\tr})$ is a non-torsion point, then
\begin{equation}\label{trinf}
\hat{h}(P)\geq \frac{1}{216\,h^*(j_E)^5},
\end{equation}
where $h^*(j_E)=h(j_E)+10$.
\end{thm}

Before we begin the proof of this result, let us first recall a few facts about elliptic curves over $\R$.  Given an elliptic curve $E/\R$, 
there exists a unique element $\tau$ in the set
\begin{equation}\label{realtau}
{\mathcal T} = \Big\{it\,\mid\,t\geq1\Big\}\cup\Big\{\frac{1}{2}+it\,\mid\,t>\frac{1}{2}\Big\},
\end{equation}
such that $j(\tau)=j_E$.  Letting $q=e^{2\pi i\tau}\in\R$, we have isomorphisms
\begin{equation}\label{compisoms}
\C^*/q^\Z \stackrel{\sim}{\longrightarrow} \C/L  \stackrel{\sim}{\longrightarrow}  E(\C);
\end{equation}
the first map in $(\ref{compisoms})$ is the inverse of the exponential map $z\mapsto e^{2\pi i z}$; while the second map in $(\ref{compisoms})$ 
can be given explicitly in terms of Weierstrass functions.  These maps restrict to isomorphisms
\begin{equation}\label{realisoms}
\R^*/q^\Z \stackrel{\sim}{\longrightarrow} {\mathcal R} \stackrel{\sim}{\longrightarrow}  E(\R),
\end{equation}
where 
\begin{equation*}
{\mathcal R}= \Big\{z\in\C/L\,\mid\,\Re(z)\in\frac{1}{2}\Z \Big\}.
\end{equation*}
(Note that, while $\Re(z)$ is not well defined on $\C/L$, the property $\Re(z)\in\frac{1}{2}\Z$ is well defined, since $\Re(\tau)\in\frac{1}{2}\Z$).  
In particular, $E(\R)$ has either one or two connected components according to whether $\Re(\tau)=1/2$ or $\Re(\tau)=0$, respectively; 
or equivalently whether $q<0$ or $q>0$, respectively.  For more details on these assertions, see \cite{SilvermanII}, \S V.2.  

\begin{lem}\label{reallemma}
Let $E/\R$ be an elliptic curve, and let $\tau=a+bi\in{\mathcal T}$ be chosen so that $j(\tau)=j_E$.  If $Z=\{P_1,\dots,P_N\}\subset E(\R)$ is a set of 
$N$ distinct real points, then we have the lower bound
\begin{equation*}
D(Z) \geq \frac{b}{4\pi|\tau|^2}e^{-8\pi|\tau|^2/bN}
\end{equation*}
for the complex discrepancy of the set $Z$.
\end{lem}
\begin{proof}
Let $\{z_1,\dots,z_N\}\subset\C$ be a set of coset representatives for the pullback of the set $Z$ under the second map in $(\ref{compisoms})$.  
Thus $z_j=r_{1,j}+r_{2,j}\tau$, where $r_{1,j}\in\frac{1}{2}\Z$ for all $j=1,\dots,N$, by $(\ref{realisoms})$.   Let $\delta$ denote the probability measure on $\C/L$ 
assigning a mass of $1/N$ at each point $z_j$.  Recall that the dual group $\Gamma_E$ is parametrized by the lattice $L$ under the correspondence $(\ref{chardef})$.  
In particular, if $\omega=n_1+n_2\tau\in L$ is a lattice point with $n_1$ even and $n_2=0$, then 
\begin{equation*}
\begin{split}
\widehat{\delta}(\gamma_\omega) & = \frac{1}{N}\sum_{j=1}^{N}\overline{\gamma_\omega(z_j)} \\
        & = \frac{1}{N}\sum_{j=1}^{N}e^{-2\pi i n_1r_{1,j}} \\
        & = 1.
\end{split}
\end{equation*}
We apply this for $\omega=\pm2+0\tau$, in which case $\omega'=\mp2\tau$.  Since $\widehat{\lambda_t}(\gamma_\omega)=\frac{b}{2\pi|\omega'|^2}e^{-2\pi t|\omega'|^2/b}$ by $(\ref{explicitcoef})$ and $(\ref{lambdat})$, we then have
\begin{equation*}
\begin{split}
D(Z) & = \sum_{\gamma\in\Gamma_E}\widehat{\lambda_{\frac{1}{N}}}(\gamma)|\widehat{\delta}(\gamma)|^2 \\
        & \geq \widehat{\lambda_{\frac{1}{N}}}(\gamma_{(2+0\tau)}) + \widehat{\lambda_{\frac{1}{N}}}(\gamma_{(-2+0\tau)})\\
        & = \frac{b}{4\pi|\tau|^2}e^{-8\pi|\tau|^2/bN},
\end{split}
\end{equation*}
which is the desired inequality.
\end{proof}

\begin{proof}[Proof of Theorem~\ref{trtheorem}]
Let us write $h=h(j_E)$ and $h^*=h+10$.  In view of Lemma~\ref{smalllem}, it suffices to show that
\begin{equation}\label{trsmallbound}
\Big|\Big\{P\in E(\Q_{\tr})\,\mid\,\h(P)\leq \frac{1}{24 h^*}\Big\}\Big| \leq 3(h^*)^2.
\end{equation}
Let $k\subset \Q_{\tr}$ be a number field over which $E$ is defined, and let $Z=\{P_1,\dots, P_N\}\subset E(\Q_{\tr})$ 
be a set of $N$ totally real points with $\h(P_j)\leq 1/24h^*$ for all $1\leq j\leq N$.  
Select a number field $K$, containing $k(P_1,\dots,P_N)$, such that $E/K$ is semistable.  Denote by $d=[K:\Q]$ and $d_v=[K_v:\Q_v]$ the global degree and local degrees of $K$, respectively.

Let $v$ be an archimedean place of $K$, corresponding to an embedding $\sigma:K\hookrightarrow\C$.  Since $k$ is totally real, we have $\sigma(k)\subset\R$, and therefore we can view $E$ as an elliptic curve defined over $\R$, with $Z\subset E(\R)$.  Select $\tau_v=a_v+b_vi$ in the set ${\mathcal T}$ defined in $(\ref{realtau})$, 
with $\sigma(j_E)=j(\tau_v)$.  By Lemma~\ref{reallemma}, we have 
\begin{equation}\label{locallowerbound}
D_v(Z)\geq \frac{1}{N}\phi\Big(\frac{8\pi|\tau_v|^2}{b_vN}\Big),
\end{equation}
where $\phi(x) = \frac{2e^{-x}}{x}$.

In order to assemble these local estimates at the archimedean places we will apply Jensen's inequality
\begin{equation}\label{ji}
\sum_v w_v\phi(x_v) \geq \phi\Big(\sum_v w_v x_v\Big)
\end{equation}
to the convex function $\phi(x)$; here $\{w_v\}$ is a finite set of positive weights with $\sum_v w_v=1$, and the $x_v$ are arbitrary positive numbers.  Applying $(\ref{ji})$ with $w_v=d_v/d$ and $x_v=8\pi|\tau_v|^2/b_vN$, and using the local lower bound $(\ref{locallowerbound})$, we have
\begin{equation}\label{jensenapp}
\begin{split}
\D(Z) & \geq \sum_{v\in\M_K^{\infty}}\frac{d_v}{d}D_v(Z) \\
        & \geq \frac{1}{N}\sum_{v\in\M_K^{\infty}}\frac{d_v}{d}\phi\Big(\frac{8\pi|\tau_v|^2}{b_vN}\Big) \\
        & \geq \frac{1}{N}\phi\Big(\sum_{v\in\M_K^{\infty}}\frac{d_v}{d}\cdot\frac{8\pi|\tau_v|^2}{b_vN}\Big) \\
        & \geq \frac{1}{N}\phi(t_N),
\end{split}
\end{equation}
where $t_N=4h^*/N$.  The last inequality in $(\ref{jensenapp})$ requires an explanation: first note that since $a_v \in \{ 0, 1/2 \}$ and $b_v\geq1/2$ for all archimedean places $v$ of $K$, we have
\begin{equation*}
\begin{split}
8\pi|\tau_v|^2/b_v & = 8\pi(a_v^2+b_v^2)/b_v \\
        & \leq 8\pi b_v + 4\pi.
\end{split}
\end{equation*}
Therefore, by Lemma \ref{jinvlem} of Appendix \ref{appendix}, we have
\begin{equation*}
\begin{split}
\sum_{v\in\M_K^{\infty}}\frac{d_v}{d}\cdot\frac{8\pi|\tau_v|^2}{b_vN}
        & \leq \frac{1}{N}\sum_{v\in\M_K^{\infty}}\frac{d_v}{d}(8\pi b_v + 4\pi ) \\
        & \leq \frac{1}{N}\sum_{v\in\M_K^{\infty}}\frac{d_v}{d}(4\log^+|j_E|_v+24+4\pi) \\
        & \leq (4h+24+4\pi)/N \\
        & < 4h^*/N,
\end{split}
\end{equation*}
which, in combination with the fact that $\phi(x)$ is decreasing for $x>0$, establishes the last inequality in $(\ref{jensenapp})$.

Now, combining the lower bound $(\ref{jensenapp})$ on $\D(Z)$ with the upper bound on $\D(Z)$ given by Theorem \ref{mainthm}, we have
\begin{equation}\label{ftbound}
\frac{2}{t_N}e^{-t_N} \leq \frac{N}{6h^*} + \frac{1}{2}\log N + \frac{1}{12}h + \frac{16}{5},
\end{equation}
since $\h(Z)\leq 1/24h^*$.  We will use this inequality to show that
\begin{equation}\label{abbound}
N \leq A\log N + B,
\end{equation}
where $A=2h^*$, and $B=4h^*(\frac{1}{12}h + \frac{16}{5})$.  First, note that $(\ref{abbound})$ is trivial if $N\leq B$, so we may assume that $N>B\geq4\cdot10\cdot\frac{16}{5}=128$ (note that $h^*\geq10$).  If $t_N>\log(6/5)$, then
\begin{equation*}
N < \frac{4h^*}{\log(6/5)} < (2\log(128)+64/5)h^* < A\log N + B.
\end{equation*}
On the other hand, if $t_N\leq \log(6/5)$, then by $(\ref{ftbound})$ we have 
\begin{equation*}
\begin{split}
\frac{5N}{12h^*} & = \frac{5}{3t_N} \\
        & \leq \frac{2}{t_N}e^{-t_N} \\
        & \leq \frac{N}{6h^*} + \frac{1}{2}\log N + \frac{1}{12}h + \frac{16}{5}.
\end{split}
\end{equation*}
Therefore
\begin{equation*}
\frac{N}{4h^*} \leq \frac{1}{2}\log N + \frac{1}{12}h + \frac{16}{5}
\end{equation*}
and $(\ref{abbound})$ follows.  

Finally, combining $(\ref{abbound})$, Lemma~\ref{tlem}, and some elementary calculus, we have
\begin{equation*}
\begin{split}
N & \leq \frac{e}{e-1}\Big(2h^*\log(2h^*) + 4h^*(\frac{1}{12}h + \frac{16}{5})\Big) \\
        & = \frac{e}{e-1}\Big(2h^*\log(2h^*) + 4h^*(\frac{1}{12}h^* + \frac{142}{60})\Big) \\
        & = \frac{e}{e-1}\Big(\frac{2\log(2h^*)}{h^*} + \frac{1}{3} + \frac{568}{60h^*}\Big)(h^*)^2 \\
        & \leq \frac{e}{e-1}\Big(\frac{2\log(2(10))}{10} + \frac{1}{3} + \frac{568}{60(10)}\Big)(h^*)^2 \\
        & < 3(h^*)^2,
\end{split}
\end{equation*}
since $h^*\geq10$.  This establishes $(\ref{trsmallbound})$, and completes the proof of the theorem. 
\end{proof}

\end{subsection}

\begin{subsection}{Small Cyclotomic Points}\label{app2}
In \cite{Zhang}, S.~Zhang notes that the Szipro-Ullmo-Zhang equidistribution theorem \cite{SzpiroUllmoZhang}, together with a
restriction of scalars argument, yields the following result: If $k$ is a totally real number field and $A / k$ is an abelian variety, then
there exists an $\varepsilon > 0$ such that $\hhat(P) \geq \varepsilon$ for all but finitely many points $P \in A(k(\mu_\infty))$, where
$k(\mu_\infty)$ denotes the maximal cyclotomic extension of $k$.
The key point in the proof is that the maximal totally real subfield $k(\mu_{\infty})^+$ of $k(\mu_{\infty})$ has index 2.  In this section, we prove an explicit quantitative version of this result when $A=E$ is an elliptic curve.

\begin{thm}
\label{Qabtheorem}
Let $k$ be a totally real subfield of $\Qbar$, and let $k(\mu_{\infty})$ denote the maximal cyclotomic extension of $k$.
If $E/k$ is an elliptic curve with j-invariant $j_E$, then 
\begin{equation}\label{cycltor}
|E(k(\mu_\infty))_{tor}| \leq 36\,h^*(j_E)^4,
\end{equation}
\begin{equation}\label{cyclliminf}
\liminf\{\h(P)\,\mid\,P\in E(k(\mu_\infty))\} \geq \frac{1}{96\,h^*(j_E)},
\end{equation}
and if $P\in E(k(\mu_\infty))$ is a non-torsion point, then 
\begin{equation}\label{cyclinf}
\h(P)\geq \frac{1}{864\,h^*(j_E)^5},
\end{equation}
where $h^*(j_E)=h(j_E)+10$.
\end{thm}

\begin{proof}
Suppose $E/k$ is given in Weierstrass normal form by the equation $y^2 = f(x)$.
If $\tau \in \Gal(k(\mu_\infty) / k)$ is any complex conjugation automorphism, then
$(k(\mu_\infty))^\tau = k(\mu_{\infty})^+$ is a totally real field.  (This follows from
the fact that $\Gal(k(\mu_\infty) / k)$ is abelian.)

Let $P \in E(k(\mu_\infty))$, and define $P_1 = P + P^\tau$, $P_2 = P - P^\tau$.  
Then $P_1$ is totally real (i.e., $P_1^\tau = P_1$) and
$P_2$ is totally imaginary (i.e., $P_2^\tau = -P_2$).

Note that the totally imaginary points on $E$ are in bijection with the totally real points
on the quadratic twist $E'$ of $E$ defined by the Weierstrass equation
\[
E' : y^2 = -f(x).
\]
Note also that $j_{E'} = j_E$, since $E$ and $E'$ are isomorphic over $k(\sqrt{-1})$.

In order to prove $(\ref{cycltor})$, note that by Theorem~\ref{trtheorem}, there are at most $M = 3 h^*(j_E)^2$ totally real torsion points on $E$.  
Applying the same result to $E'$, and using the fact that $j_{E'} = j_E$, we see
that there are at most $M$ totally imaginary torsion points on $E$ as well.
If $P \in E(k(\mu_\infty))_{\tor}$, then as $2P = P_1 + P_2$, 
it follows that there are at most $M^2$ possibilities for $2P$, and
thus at most $4M^2$ possibilities for $P$.  This proves the bound $(\ref{cycltor})$.

We now turn to the proof of $(\ref{cyclliminf})$ and $(\ref{cyclinf})$.  
First, we claim that given a non-torsion point $P\in E(k(\mu_\infty))$, either there exists a non-torsion totally real point $Q$ on $E$ with $\h(P)\geq\h(Q)/4$, 
or there exists a non-torsion totally real point $Q'$ on $E'$ with $\h(P)\geq\h(Q')/4$.  
To see this, note that by the parallelogram law we have $\h(P_1)+\h(P_2)=2\h(P)+2\h(P^\tau)=4\h(P)>0$.  
Therefore either $\h(P_1)>0$, in which case we take $Q=P_1$; or else $\h(P_2)>0$, in which case we take $Q'=P_2'$, 
the image of $P_2$ under the isomorphism $E\simeq E'$.

In view of this claim, $(\ref{cyclliminf})$ follows immediately from $(\ref{trliminf})$, and $(\ref{cyclinf})$ from $(\ref{trinf})$.
\end{proof}

\end{subsection}

\begin{subsection}{Small Totally $p$-adic Points}\label{app3}

Let $p$ be a prime number.  A subfield $L$ of $\Qbar$ is called {\em totally $p$-adic} if the rational prime 
$p$ splits completely in every number field $L' \subseteq L$, 
or equivalently, if every embedding of $L$ into $\overline{\QQ}_p$ has image contained in $\QQ_p$.

The following result is an analogue for elliptic curves of a result of Bombieri and Zannier \cite{BombieriZannier} for the multiplicative group:

\begin{thm}\label{tpadictheorem}
Let $k$ be a number field, let $p\neq 2$ be a prime number, and let $E/k$ be an elliptic curve
having semistable reduction at all places of $k$ lying over $p$.  If $L/k$ is totally $p$-adic algebraic extension, then
\begin{equation}\label{tpadictor}
|E(L)_{\tor}| \leq \frac{8M}{5 \log p} \Big( \log M + 2\log p + h^\dagger(j_E) \Big),
\end{equation}
\begin{equation}\label{tpadicliminf}
\liminf\{\h(P)\,\mid\,P\in E(L)\} \geq \frac{\log p}{8M},
\end{equation}
and if $P\in E(L)$ is a non-torsion point, then
\begin{equation}\label{tpadicinf}
\hat{h}(P)\geq \frac{25}{512}\Big(\frac{\log p}{M}\Big)^3\Big( 
\log M + 2\log p + h^\dagger(j_E) \Big)^{-2},
\end{equation}
where 
\[
M = \max \{ p + 1 + 2\sqrt{p}, 12 \nu \}, \; h^\dagger(j_E) = \frac{1}{6} h(j_E) + \frac{32}{5},
\]
and $\nu$ is the maximum over all places $w \in M_k$ lying over $p$ of the quantity 
$w^+(j_E)=\max\{0,-{\rm ord}_w(j_E)\}$.
\end{thm}

\begin{proof}
In view of Lemma~\ref{smalllem}, it suffices to show that
\begin{equation}\label{tpadicsmall}
\Big|\Big\{P\in E(L)\,\mid\,\h(P)\leq\frac{\log p}{8M}\Big\}\Big| \leq \frac{8M}{5 \log p} \Big( 
\log M + 2\log p + \frac{1}{6} h(j_E) + \frac{32}{5}
\Big).
\end{equation}

Let $Z = \{ P_1,\ldots,P_N \}$ be a distinct set of points of $E(L)$, with $\h(P_j)\leq\frac{\log p}{8M}$ for all $1\leq j\leq N$, and choose a finite extension $K$ of $k$ over which all points of $Z$ are defined.
Let $v$ be a place of $K$ lying over $p$, and note that if $\beta \in K$ is a nonzero element with $|\beta|_v < 1$, 
then in fact $|\beta|_v \leq p^{-1}$.  

\medskip

Suppose first that $E$ has good reduction at $v$, and let $\bar{E}$ be its reduction, which is
an elliptic curve over the residue field $\FF_v$ of $\O_K$ at $v$.
Setting $r = p^{-1}$, we recall from (\ref{epdisc}) that
\[
D_i(Z) \geq (D_{i,r}(Z) + \frac{1}{N})\log(1/r) .
\]

Let $m = |\bar{E}(\FF_v)|$, which by Hasse's theorem satisfies the inequality
$m \leq p + 1 + 2\sqrt{p}$.  Then by the pigeonhole principle, we have
\begin{equation*}
\begin{split}
D_{i,r}(Z) &= \sum_{\alpha\in\bar{E}(\F_v)}\bigg(\frac{N_r(\alpha)}{N}\bigg)^2 \\
&\geq \sum_{\alpha\in\bar{E}({\F}_v)} \bigg(\frac{N/m}{N}\bigg)^2 \\
& = \frac{1}{m}.
\end{split}
\end{equation*}
Thus 
\[
D_i(Z) \geq \left( \frac{1}{m} - \frac{1}{N} \right)\log p .
\]

\medskip

Now suppose that $E$ has multiplicative reduction at $v$, and let $\nu_v = -\ord_v(j_E)\geq1$.
Recall that the retraction homomorphism $r:E(K_v)\to\R/\Z$ factors as
\begin{equation*}
E(K_v)\to K_v^*/q^\Z\to\R/\Z,
\end{equation*}
where $q\in K_v^*$ with $|q|_v=|1/j_E|_v$; here the first map is the Tate parametrization, 
and the second map is given by $u\mapsto\log|u|_v/\log|q|_v$ (cf. $\S$\ref{RetractionHomomorphismSection}).  
Note that $\text{Im}(r)=\langle1/\nu_v\rangle\subset\R/\Z$, and therefore $e^{2\pi inr(P_j)}=1$ whenever $\nu_v\mid n$.   
Letting $D_j(Z)$ denote the local retraction discrepancy as defined in Section~\ref{discnonarch}, it follows 
from Proposition~\ref{retprop} that
\begin{equation}
\label{RetractionDiscrepancyEstimates}
\begin{split}
D_j(Z) & \geq \frac{\log|j_E|_v}{4\pi^2}\sum_{n\in\Z\setminus\{0\}}\frac{1}{n^2}\Big|\frac{1}{N}\sum_{1\leq j\leq N}e^{2\pi inr(P_j)}\Big|^2 \\
        & \geq \frac{\log|j_E|_v}{4\pi^2}\sum_{n\in \nu_v\Z\setminus\{0\}}\frac{1}{n^2} \\
        & =\frac{\log|j_E|_v}{12\nu_v^2} \\
        & = \frac{\log p}{12\nu_v}.
\end{split}
\end{equation}

\medskip

Using the definition of the global discrepancy and of the constant $M$, together with the fact that 
\[
\sum_{\stackrel{v \in M_K}{v\mid p}} \frac{d_v}{d} = 1,
\]
we see in all cases that the global discrepancy $\D(Z)$ of $Z$ satisfies
\[
\D(Z) \geq \Big( \frac{1}{M} - \frac{1}{N} \Big)\log p.
\]
Combining this with (\ref{mainbound}), we find that
\begin{equation}\label{tpadicmainbound}
\begin{split}
\Big( \frac{1}{M} - \frac{1}{N} \Big)\log p
        & \leq 4\h(Z) + \frac{1}{N}\Big(\frac{1}{2}\log N + \frac{1}{12}h(j_E) + \frac{16}{5}\Big) \\
        & \leq \frac{\log p}{2M} + \frac{1}{N}\Big(\frac{1}{2}\log N + \frac{1}{12}h(j_E) + \frac{16}{5}\Big),
\end{split}
\end{equation}
since $\h(Z)\leq\frac{\log p}{8M}$.  It follows that
\begin{equation}
\label{tpadicmainbound2}
N \leq \frac{M}{\log p} \Big( \log N + 2\log p + \frac{1}{6}h(j_E) + \frac{32}{5} \Big).
\end{equation}
Applying Lemma~\ref{tlem} with $A = \frac{M}{\log p}$ and 
$B = \frac{M}{\log p} \Big( 2\log p + \frac{1}{6}h(j_E) + \frac{32}{5} \Big)$,
we obtain the upper bound $N \leq \frac{8}{5} (A\log A + B)$.  
The inequality $(\ref{tpadicsmall})$ follows upon noting that $\log\log p \geq 0$.
\end{proof}

By suitably modifying the above argument, one can establish the following generalization of Theorem~\ref{tpadictheorem}; we omit the details of the proof.  
For the statement, we say a subfield $L$ of $\Qbar$ is {\em totally $p$-adic of type $(\ee,\ff)$} if for any embedding of $L$ into $\overline{\QQ}_p$, 
the image of $L$ is contained in a finite extension of $\QQ_p$ whose absolute ramification and residue degrees are bounded by $\ee$ and $\ff$, respectively.  

\begin{thm}\label{tpadictheoremalternate}
Let $k$ be a number field, let $p\neq 2$ be a prime number, and let $E/k$ be an elliptic curve
having semistable reduction at all places of $k$ lying over $p$.
If $L/k$ is an algebraic extension which is totally $p$-adic of type $(e,f)$ and $q = p^f$, then
\begin{equation}
|E(L)_{tor}| \leq \frac{8 e M}{5 \log p} \Big(\log(eM) + \frac{2\log p}{e} + h^\dagger(j_E) \Big),
\end{equation}
\begin{equation}
\liminf\{\h(P)\,\mid\,P\in E(L)\} \geq \frac{\log p}{8eM},
\end{equation}
and if $P\in E(L)$ is a non-torsion point, then
\begin{equation}
\hat{h}(P)\geq \frac{25}{512}\Big(\frac{\log p}{eM}\Big)^3\Big( 
\log(eM) + \frac{2\log p}{e} + h^\dagger(j_E) \Big)^{-2},
\end{equation}
where 
\[
M = \max \{ q + 1 + 2\sqrt{q}, 12e \nu \}, \; h^\dagger(j_E) = \frac{1}{6} h(j_E) + \frac{32}{5},
\]
and $\nu$ is the maximum over all places $w \in M_k$ lying over $p$ of the quantity
$w^+(j_E)=\max\{0,-\ord_w(j_E)\}$.
\end{thm}

\end{subsection}

\end{section}


\appendix

\begin{section}{A quantitative refinement of Elkies' theorem}\label{appendix}

Our goal in this appendix is to give a proof of Proposition~\ref{elkies}.  
The proof is simpler than the one given in \cite{LangAT}, $\S$VI, Theorem 5.1, in that
we use only basic results from Fourier analysis, whereas the proof presented by Lang
uses a number of nontrivial results concerning elliptic differential equations and 
eigenfunctions of the Laplacian on a Riemannian manifold.
At the same time, our proof yields better quantitative information than one obtains in the
general case by using certain explicit estimates for the modular $j$-function.  We also include a 
discrepancy term in the estimate which plays an important role in this paper.

\medskip

For the reader's convenience, we recall the statement of Proposition~\ref{elkies}. 

\medskip

{\bf Proposition~\ref{elkies}.} 
Let $E/\C$ be an elliptic curve with j-invariant $j_E$, and let 
\begin{equation*}
Z=\{P_1,\dots,P_N\}\subset E(\C)
\end{equation*}
be a set of $N$ distinct points.  Then $D(Z)>0$, and 
\begin{equation*}
\sum_{\stackrel{1\leq i,j\leq N}{i\neq j}}\lambda(P_i-P_j)\geq N^2D(Z) - \frac{N\log N}{2} -\frac{N}{12}\log^+|j_E| - \frac{16}{5}N.
\end{equation*}

\medskip

Before giving the proof, we need a series of preliminary results.  

\begin{lem}
The kernel $g_t$ is nonnegative.
\end{lem}
\begin{proof}
Let $\psi\in \Sm(E(\C))$ be a nonnegative test function, normalized so that $\widehat{\psi}(\gamma_0)=\int_{E(\C)}\psi(P)d\mu(P)=1$, 
and define $u(P,t)\in C^\infty\big(E(\C)\times[0,+\infty)\big)$ by 
\begin{equation*}
\begin{split}
u(P,t) & = 
\begin{cases}
g_t*\psi(P) & \text{ if } t>0 \\
\psi(P) & \text{ if } t=0.
\end{cases} \\
        & = \sum_{\gamma\in\Gamma_E}\widehat{\psi}(\gamma)e^{-t/\widehat{\lambda}(\gamma)}  \gamma(P).
\end{split}
\end{equation*}
In view of Proposition \ref{eigenfunction} we see that $u$ is a solution to the heat equation
\begin{equation*}
\frac{\partial u}{\partial t} = \Delta u.
\end{equation*}
Such a function satisfies a ``maximum principle,'' meaning it must take its extrema on the boundary $E(\C)\times\{0\}$ (cf. \cite{Folland}, (4.16)).  
In particular, we have
\begin{equation*}
\inf_{P,t}u(P,t) = \inf_{P}u(P,0),
\end{equation*}
and since $u(P,0)=\psi(P)\geq0$, it follows that $u(P,t)$ is nonnegative.  
If we now fix $t>0$ and take a sequence of such $\psi$ so that $\psi(P)d\mu(P)$ converges weakly to a unit point mass at $P=O$, 
we have $u(P,t)\to g_t(P)$, and therefore $g_t(P)\geq0$.  
\end{proof}

\begin{lem}\label{pl1}
For $P\in E(\C)\setminus\{O\}$ and $t>0$, we have
\begin{equation*}
\lambda(P)\geq\lambda_t(P)-t.
\end{equation*}
\end{lem}
\begin{proof}
We have
\begin{equation*}
\begin{split}
\lambda_{t}(P)-\lambda(P) & = \int_0^{t}\frac{\partial}{\partial s}\lambda_s(P)ds \\
        & = \int_0^{t}\frac{\partial}{\partial s}\Big(\sum_{\gamma\in\Gamma_E}\widehat{\lambda}(\gamma)\widehat{g_s}(\gamma)\gamma(P)\Big)ds \\
        & = \int_0^{t}\frac{\partial}{\partial s}\Big(\sum_{\gamma\in\Gamma_E\setminus\{\gamma_0\}}\widehat{\lambda}(\gamma)e^{-s/\widehat{\lambda}(\gamma)}  
            \gamma(P)\Big)ds \\
        & = -\int_0^{t}\Big(\sum_{\gamma\in\Gamma_E\setminus\{\gamma_0\}}e^{-s/\widehat{\lambda}(\gamma)} \gamma(P)\Big)ds \\
        & = -\int_0^{t}g_s(P)ds + t \\ 
        & \leq t,
\end{split}
\end{equation*}
since $g_s(P)\geq0$.
\end{proof}

We will also require the following estimate involving the modular j-function $j:\H\to\C$.

\begin{lem}\label{jinvlem}
For $\tau=a+bi\in\H$, we have 
\begin{equation}\label{jinv}
2\pi b \leq \log^+|j(\tau)| + 6.
\end{equation}
\end{lem}
\begin{proof}
Let $0<t<1$ be a parameter to be chosen later, and put $q=e^{2\pi i\tau}$.  If $|q|\geq t$, then 
\begin{equation}\label{nb1}
\begin{split}
2\pi b & = -\log|q| \\
        & \leq \log(1/t) \\
        & \leq \log^+|j(\tau)| + \log(1/t).
\end{split}
\end{equation}
On the other hand, suppose that $|q|<t$.  Recall that $j(\tau)=1728g_2(\tau)^3/\Delta(\tau)$, 
where $\Delta(\tau)=(2\pi)^{12}q\prod_{n\geq1}(1-q^n)^{24}$ is the modular discriminant function, and
\begin{equation*}
g_2(\tau)=\frac{(2\pi)^4}{12}\Big(1+240\sum_{d\geq1}\sum_{m\geq1}d^3q^{md}\Big)
\end{equation*}
is the Eisenstien series of weight 4 (cf. \cite{Serre}, ch. VII).  We have
\begin{equation*}
\begin{split}
\log|\Delta(\tau)/q| & = 12\log(2\pi)+ 24\sum_{n\geq1}\log|1-q^n| \\
        & = 12\log(2\pi) -24\Re\sum_{n\geq1}\sum_{k\geq1}k^{-1}q^{kn} \\
        & \leq 12\log(2\pi) +24\sum_{n\geq1}\sum_{k\geq1}k^{-1}t^{kn} \\
        & = 12\log(2\pi) +24\sum_{k\geq1}\frac{k^{-1}t^k}{1-t^k} \\
        & \leq 12\log(2\pi) +\frac{24}{1-t}\sum_{k\geq1}k^{-1}t^k \\
        & = 12\log(2\pi) +\frac{24}{1-t}\log\Big(\frac{1}{1-t}\Big).
\end{split}
\end{equation*}
Also, using the inequality $|1+u|\geq1-|u|$ and the identity $\sum_{d\geq1}d^3t^d= t(1+4t+t^2)/(1-t)^4$, we estimate
\begin{equation*}
\begin{split}
|g_2(\tau)| & = \frac{(2\pi)^4}{12}\Big|1+240\sum_{d\geq1}\sum_{m\geq1}d^3q^{md}\Big| \\
        & \geq \frac{(2\pi)^4}{12}\Big(1-240\sum_{d\geq1}\sum_{m\geq1}d^3t^{md}\Big) \\
        & = \frac{(2\pi)^4}{12}\Big(1-240\sum_{d\geq1}d^3\frac{t^d}{1-t^d}\Big) \\
        & \geq \frac{(2\pi)^4}{12}\Big(1-\frac{240}{1-t}\sum_{d\geq1}d^3t^d\Big) \\
        & = \frac{(2\pi)^4}{12}\Big(1-\frac{240t(1+4t+t^2)}{(1-t)^5}\Big). \\
\end{split}
\end{equation*}
Combining these two estimates, we have
\begin{equation}\label{nb2}
\begin{split}
2\pi b & = -\log|q| \\
        & = \log|j(\tau)| - \log|1728g_2(\tau)^3q/\Delta(\tau)| \\
        & \leq \log^+|j(\tau)| + \frac{24}{1-t}\log\Big(\frac{1}{1-t}\Big) - \log\Big(1-\frac{240t(1+4t+t^2)}{(1-t)^5}\Big).
\end{split}
\end{equation}

We now select $t$ to be the (only) root $t_0$ of the polynomial $f(t)=240t(1+4t+t^2)-(1-t)^6$ in the interval $(0,1)$.  This choice, while perhaps not quite optimal, does simplify matters since it satisfies the identity
\begin{equation}\label{tzero}
\log(1/t_0) = -\log\Big(1-\frac{240t_0(1+4t_0+t_0^2)}{(1-t_0)^5}\Big).
\end{equation}
Therefore, in view of the two cases $(\ref{nb1})$ and $(\ref{nb2})$, and the identity $(\ref{tzero})$, it follows that
\begin{equation}\label{c1}
2\pi b \leq \log^+|j(\tau)| + \log(1/t_0) + \frac{24}{1-t_0}\log\Big(\frac{1}{1-t_0}\Big)
\end{equation}
for all $\tau=a+ib\in\H$.  Finally, using the fact that $f(t)$ is increasing for $0< t< 1$, it is straightforward to show that $1/250\leq t_0\leq1/249$, and therefore 
\begin{equation*}
\begin{split}
\log(1/t_0) + \frac{24}{1-t_0}\log\Big(\frac{1}{1-t_0}\Big) & \leq \log(250) + \frac{24}{1-1/249}\log\Big(\frac{1}{1-1/249}\Big) \\
        & < 6,
\end{split}
\end{equation*}
which, along with $(\ref{c1})$, finishes the proof of the lemma.
\end{proof}

\begin{remark}
In \cite{HindrySilvermanI}, Hindry and Silverman prove the reverse inequality
$2\pi b \geq \log^+|j(\tau)| - 2.304$
by a similar method.
\end{remark}

\begin{lem}\label{pl2}
Let $E/\C$ be an elliptic curve with j-invariant $j_E$.  Then for $0<t\leq1$, we have
\begin{equation}\label{sh}
\lambda_t(O)\leq  \frac{1}{2}\log(1/t) + \frac{1}{12}\log^+|j_E| + 11/5.
\end{equation}
\end{lem}
\begin{proof}
Fix a parametrization $E(\C)\simeq\C/L$, where $L=\Z+\tau\Z$ is a lattice and $\tau=a+bi$ is chosen in the usual fundamental domain 
\begin{equation*}
{\mathcal F}=\big\{z\,\mid\,-1/2<\Re(z)\leq1/2,|z|>1\big\}\cup\big\{e^{i\theta}\,\mid\,\pi/3\leq\theta\leq\pi/2\big\}
\end{equation*}
for the action of the modular group on $\H$.  Thus $j_E=j(\tau)$, and note in particular that $b\geq\sqrt{3}/2$.

Given a nonzero lattice point $\omega\in L$, we have the associated (nontrivial) character $\gamma_\omega\in\Gamma_E$, as defined in $(\ref{chardef})$.  
It follows from $(\ref{lapchar})$ and Proposition \ref{eigenfunction} that
\begin{equation*}
\widehat{\lambda}(\gamma_\omega)=\frac{b}{2\pi|\omega'|^2},
\end{equation*}
where $\omega'$ is defined in $(\ref{latperm})$.  

Since $\omega\mapsto\omega'$ is a permutation on $L$ which fixes zero, we have
\begin{equation*}
\begin{split}
\lambda_{t}(O) & = \sum_{\gamma\in\Gamma_E\setminus\{\gamma_0\}}\widehat{\lambda}(\gamma)e^{-t/\widehat{\lambda}(\gamma)} \\
        & = \sum_{\omega\in L\setminus(0)}\frac{b}{2\pi|\omega|^2}e^{-2\pi|\omega|^2 t/b} \\
        & \leq \frac{b}{2\pi}\sum_{k\in\Z\setminus(0)}\frac{1}{k^2} +  \frac{b}{\pi}\sum_{\omega\in L\cap\H}\frac{1}{|\omega|^2}e^{-2\pi|\omega|^2 t/b} \\
        & = \frac{\pi b}{6} + \frac{b}{\pi}\int_{0}^{\infty}\frac{1}{x}e^{-2\pi xt/b}dS(x),
\end{split}
\end{equation*}
where the latter is a Riemann-Stieltjes integral, and
\begin{equation*}
S(x) = \sum_{\stackrel{\omega\in L\cap\H}{|\omega|\leq \sqrt{x}}}1
\end{equation*}
is the number of lattice points $\omega\in L$, with positive imaginary part, in the disc $\{|z|\leq\sqrt{x}\}$.  
An elementary lattice point counting argument provides a bound $S(x)\leq \frac{\pi}{2b}x+\frac{1}{b}\sqrt{x}$.  
Also, it is clear that $S(x)=0$ for $x<b^2$, since no lattice point in $L$ can have a positive imaginary part smaller than $b$.  
Integrating by parts, we have
\begin{equation}\label{b1}
\begin{split}
\lambda_{t}(O) & \leq \frac{\pi b}{6} + \frac{b}{\pi}\int_{0}^{\infty}\frac{1}{x}e^{-2\pi xt/b}dS(x)  \\
        & = \frac{\pi b}{6} + \frac{b}{\pi}\int_0^\infty S(x)\Big(\frac{1+2\pi tx/b}{x^2}\Big)e^{-2\pi tx/b}dx \\
        & \leq \frac{\pi b}{6} + \frac{b}{\pi}\int_{b^2}^\infty(\frac{\pi}{2b}x+\frac{1}{b}x^{1/2})\Big(\frac{1+2\pi tx/b}{x^2}\Big)e^{-2\pi tx/b}dx \\
        & = \frac{\pi b}{6} + \int_{b^2}^\infty\Big(\frac{1+2\pi tx/b}{2x}\Big)e^{-2\pi tx/b}dx + 
            \int_{b^2}^\infty\Big(\frac{1+2\pi tx/b}{\pi x^{3/2}}\Big)e^{-2\pi tx/b}dx \\
        & = \frac{\pi b}{6} + I_1+I_2,
\end{split}
\end{equation}
where $I_1$ and $I_2$ are intergals which we now estimate.  
First, upon making the change of variables $u=tx/b^2$, we have
\begin{equation*}
\begin{split}
I_1 & =  \frac{1}{2}\int_{b^2}^\infty x^{-1}(1+2\pi tx/b)e^{-2\pi tx/b}dx \\
        & = \frac{1}{2}\int_{t}^\infty u^{-1}(1+2\pi bu)e^{-2\pi bu}du \\
        & = \frac{1}{2}\int_{t}^1 u^{-1}e^{-2\pi bu}du +  \frac{1}{2}\int_{1}^\infty u^{-1}e^{-2\pi bu}du +  \pi b\int_{t}^\infty e^{-2\pi bu}du \\
        & \leq \frac{1}{2}\int_{t}^1 u^{-1}du +  \frac{1}{2}\int_{1}^\infty e^{-2\pi bu}du +  \pi b\int_{0}^\infty e^{-2\pi bu}du \\
        & = \frac{1}{2}\log(1/t) + \frac{1}{4\pi b} +  \frac{1}{2}.
\end{split}
\end{equation*}
Also, 
\begin{equation*}
\begin{split}
I_2 & = \int_{b^2}^\infty\Big(\frac{1+2\pi tx/b}{\pi x^{3/2}}\Big)e^{-2\pi tx/b}dx \\
        & \leq \frac{1}{\pi}\int_{b^2}^\infty x^{-3/2}dx + \frac{2t}{b} \int_{b^2}^\infty x^{-1/2}e^{-2\pi tx/b}dx \\
        & \leq \frac{2}{\pi b} + \frac{2t}{b^2} \int_{b^2}^\infty e^{-2\pi tx/b}dx \\
        & \leq \frac{3}{\pi b}.
\end{split}
\end{equation*}
Combining the estimates on these integrals with Lemma \ref{jinvlem}, and using the fact that $b\geq\sqrt{3}/2$ (since $\tau\in{\mathcal F}$), we deduce that 
\begin{equation*}
\begin{split}
\lambda_{t}(O) & \leq \frac{1}{2}\log(1/t) + \frac{\pi b}{6} + \frac{13}{4\pi b} +  \frac{1}{2} \\
        & \leq \frac{1}{2}\log(1/t) + \frac{1}{12}\log^+|j(\tau)| + \frac{13}{2\pi\sqrt{3}} + 1,
\end{split}
\end{equation*}
and $(\ref{sh})$ follows, since $13/2\pi\sqrt{3} + 1<11/5$.
\end{proof}

\begin{proof}[Proof of Proposition \ref{elkies}]
The positivity of the discrepancy follows at once from $(\ref{napars})$.  For $t>0$ we have
\begin{equation}\label{elkies1}
\begin{split}
\sum_{\stackrel{1\leq i,j\leq N}{i\neq j}}\lambda(P_i-P_j) & \geq \sum_{\stackrel{1\leq i,j\leq N}{i\neq j}}\lambda_t(P_i-P_j) - N^2t \\
        & = \sum_{1\leq i,j\leq N}\lambda_t(P_i-P_j) - N\lambda_t(O) - N^2t,
\end{split}
\end{equation}
by Lemma \ref{pl1}.  Selecting $t=1/N$ and using Lemma \ref{pl2}, we deduce the desired bound
\begin{equation}\label{elkies2}
\begin{split}
\sum_{\stackrel{1\leq i,j\leq N}{i\neq j}}\lambda(P_i-P_j) & \geq N^2 D(Z) - N\lambda_{1/N}(O) - N \\
        & \geq N^2 D(Z) - \frac{N\log N}{2} - \frac{N}{12}\log^+|j_E| - \frac{16}{5}N.
\end{split}
\end{equation}
\end{proof}

\end{section}


\end{document}